\documentclass[a4paper]{article}

\usepackage{geometry}
\usepackage{amsmath,amssymb,amsthm}
\usepackage{latexsym}
\usepackage{indentfirst}
\usepackage{graphicx}
\usepackage[colorlinks=true]{hyperref}
\usepackage{mathrsfs}
\usepackage{chngcntr}
\usepackage{xcolor}

\numberwithin{equation}{section}

\newtheorem{theorem}{Theorem}[section]
\newtheorem{lemma}{Lemma}[section]
\newtheorem{proposition}{Proposition}[section]

\newtheorem{definition}{Definition}[section]

\begin{document}
	
	\title{Existence and concentration of multiple solutions for the fractional logarithmic Schr\"odinger-Poisson system via the penalization method}
	
	\author{
		{Jiao Luo, Zhipeng Yang}\thanks{Email: yangzhipeng326@163.com}\\
		\small Department of Mathematics, Yunnan Normal University, Kunming, China.\\
		\small Yunnan Key Laboratory of Modern Analytical Mathematics and Applications, Kunming, China.
	}
	
	\date{}
	\maketitle
	
	\begin{abstract}
		We study positive solutions of the fractional logarithmic Schr\"odinger-Poisson system
		\[
		\begin{cases}
			\varepsilon^{2\alpha}(-\Delta)^\alpha u+V(x)u+\phi u=u\log u^2+|u|^{q-2}u, & \text{in } \mathbb{R}^3,\\
			\varepsilon^{2\alpha}(-\Delta)^\alpha \phi=u^2, & \text{in } \mathbb{R}^3,
		\end{cases}
		\]
		where \(\varepsilon>0\), \(3/4<\alpha<1\), and \(4<q<6/(3-2\alpha)\). We assume that the potential is continuous and nonnegative and attains its global minimum inside a bounded domain, with strictly larger values on the boundary. We work in the intersection of a weighted fractional Sobolev space and an Orlicz space, so that the logarithmic energy is finite and the associated functional is of class \(C^1\). For sufficiently small \(\varepsilon\), we obtain at least as many distinct positive solutions as the relative Lusternik-Schnirelmann category of the minimum set in a small neighborhood. The global maximum points of these solutions approach that minimum set as \(\varepsilon\to0\). The proof combines penalization, an autonomous Nehari-Pohozaev argument, a local natural constraint, and category theory.
	\end{abstract}
	
	
	
	\section{Introduction}
	
	Fractional Schr\"odinger equations describe models in which the kinetic energy has a nonquadratic Fourier symbol. In Laskin's fractional quantum mechanics \cite{Laskin2000}, replacing Brownian paths by L\'evy paths leads to a fractional kinetic operator. With a normalized coefficient, this operator satisfies
	\[
	\widehat{(-\Delta)^\alpha v}(\xi)=|\xi|^{2\alpha}\widehat v(\xi), \qquad 0<\alpha<1.
	\]
	Coupling the wave function to a potential generated by its own density leads to a Schr\"odinger-Poisson system. A time-dependent model considered here is
	\[
	\begin{cases}
		i \varepsilon \dfrac{\partial \Psi}{\partial t} = \varepsilon^{2\alpha}(-\Delta)^\alpha \Psi+\widetilde{V}(x)\Psi+\phi\Psi-f(\Psi), & \text{in } \mathbb{R}^{3}\times \mathbb{R},\\
		\varepsilon^{2\alpha}(-\Delta)^\alpha \phi=|\Psi|^2, & \text{in } \mathbb{R}^{3},
	\end{cases}
	\]
	where \(\varepsilon>0\), \(\alpha\in\left(\frac{3}{4},1\right)\), \(\widetilde{V}(x)=V(x)+E\) with \(E\in\mathbb R\), and \(f(e^{i\theta}\xi)=e^{i\theta}f(\xi)\).
	
	The logarithmic response has a separate origin. Bia\l ynicki-Birula and Mycielski \cite{BBM1976} introduced a logarithmic nonlinearity in nonlinear wave mechanics. In nonlinear optics, Buljan et al.\ \cite{Buljan2003} analytically studied localized white light beams in logarithmically saturable media with a noninstantaneous response. Logarithmic terms have also been considered in effective quantum models; see \cite{MR2740900}. These works motivate the study of logarithmic responses, but they do not derive the full fractional coupled system considered below. Here the positive Poisson term describes a repulsive self-interaction, whereas the additional power term gives a local focusing contribution whose energy has order \(q>4\).
	
	If the self-consistent field coupling is omitted and \(f(\Psi)=\Psi\log |\Psi|^2+|\Psi|^{q-2}\Psi\), the corresponding scalar equation is
	\[
	i \varepsilon \frac{\partial \Psi}{\partial t} = \varepsilon^{2\alpha}(-\Delta)^\alpha \Psi+(V(x)+E)\Psi-\Psi \log |\Psi|^2-|\Psi|^{q-2}\Psi, \quad \text{in } \mathbb{R}^{3}.
	\]
	Seeking standing-wave solutions of the form
	\[
	\Psi(x,t)=\exp\left(\frac{-iEt}{\varepsilon}\right)u(x),
	\]
	one is led to the fractional elliptic equation
	\[
	\varepsilon^{2\alpha}(-\Delta)^\alpha u+V(x)u=u\log u^2+|u|^{q-2}u, \quad \text{in } \mathbb{R}^{3}.
	\]
	Existence, multiplicity, and concentration results for logarithmic Schr\"odinger equations were obtained in \cite{10.57262/die035-1112-677,MR3385171,MR4097474,MR4698809,MR4048330,MR3195154,MR3451965}.
	
	A basic fractional Schr\"odinger-Poisson model is
	\begin{equation}\label{eq1.1}
		\begin{cases}
			\varepsilon^{2s} (-\Delta )^{s}u+V(x) u+\phi u=f(u), & \text{in } \mathbb{R}^{3},\\
			\varepsilon^{2t} (-\Delta )^{t}\phi=u^2, & \text{in } \mathbb{R}^{3}.
		\end{cases}
	\end{equation}
	Existence, multiplicity, and concentration results for \eqref{eq1.1}, including problems with critical growth, have been obtained under local or global conditions on \(V\); see \cite{MR4230968,MR3716931,YangYuZhao2019,YangYuZhaoCCM2019}. Related results for the classical Schr\"odinger-Poisson system and for fractional systems under different assumptions on the potentials and nonlinearities can be found in \cite{YangYu2021,GuMuYang2023,MuYangZhang2024,MR4870288}.
	
	For further results on logarithmic nonlinearities, logarithmic interaction potentials, and fractional systems with general nonlinearities, see \cite{MR4604329,MR4782770,Teng2016,MR4844569,MR3456743}. In particular, Tao and Li \cite{MR4782770} study a logarithmic fractional system with the coupling \(-\phi u\). As they note in \cite[Remark~1.2]{MR4782770}, the negative sign is essential to their boundedness estimate for Palais-Smale sequences, which does not directly extend to the positive coupling. The present problem instead contains \(+\phi u\), which contributes a positive quartic term to the reduced energy. This sign is important in the variational argument and motivates the use of the constraints introduced in Section~\ref{Sec4}.
	
	Motivated by the above works, we consider the fractional logarithmic Schr\"odinger-Poisson system
	\begin{equation}\label{eq1.2}
		\begin{cases}
			\varepsilon^{2\alpha} (-\Delta )^{\alpha}u+V(x) u+\phi u=u \log u^{2}+f(u), & \text {in } \mathbb{R}^{3},\\
			\varepsilon^{2\alpha} (-\Delta )^{\alpha}\phi=u^2, & \text{in } \mathbb{R}^{3},
		\end{cases}
	\end{equation}
	where \(V\in C(\mathbb{R}^3,\mathbb{R})\) satisfies
	\begin{itemize}
		\item[$(V_1)$] \(\inf_{x\in\mathbb{R}^3}V(x)=V_0\ge0\);
		\item[$(V_2)$] there exists an open bounded set \(\Lambda\subset\mathbb{R}^3\) such that
		\[
		V_0=\inf_{x\in\Lambda}V(x)<\min_{x\in\partial\Lambda}V(x).
		\]
	\end{itemize}
	Without loss of generality, we assume that \(0\in\Lambda\) and \(V(0)=V_0\). Moreover, for every
	\[
	q\in\left(4,2_\alpha^*=\frac{6}{3-2\alpha}\right),
	\]
	we take
	\[
	f(u)=|u|^{q-2}u, \qquad u\in\mathbb{R}.
	\]
	
	If \((v,\varphi)\) solves \eqref{eq1.2}, set \(u(x)=v(\varepsilon x)\) and \(\phi(x)=\varphi(\varepsilon x)\). Then problem \eqref{eq1.2} is equivalent to
	\begin{equation}\label{eq1.3}
		\begin{cases}
			(-\Delta)^\alpha u+V(\varepsilon x)u+\phi u=u \log u^2+|u|^{q-2}u, & \text{in } \mathbb{R}^{3},\\
			(-\Delta)^\alpha \phi=u^2, & \text{in } \mathbb{R}^{3}.
		\end{cases}
	\end{equation}
	
	In the rescaled system \eqref{eq1.3}, the field generated by \(u^2\) is
	\[
	\phi_u^\alpha(x) =C_\alpha\int_{\mathbb R^3} \frac{u(y)^2}{|x-y|^{3-2\alpha}}\,dy, \qquad C_\alpha>0.
	\]
	Thus the field term represents a nonlocal Riesz interaction, in addition to the nonlocal kinetic energy and the local nonlinear response. The concentration result below identifies where localized standing waves form as \(\varepsilon\to0\), while the multiplicity result relates their number to the topology of the minimum set of the potential.
	
	The main difficulty in studying \eqref{eq1.3} is that the logarithmic energy need not be finite on the natural fractional Sobolev space. Indeed, the formal energy functional
	\[
	I_{\varepsilon}(u) = \frac{1}{2} \int_{\mathbb{R}^{3}}\left(|(-\Delta)^{\frac{\alpha}{2}}u|^{2}+V(\varepsilon x)|u|^{2}\right)\,dx +\frac{1}{4}\int_{\mathbb{R}^{3}}\phi_{u}^\alpha u^2\,dx -\int_{\mathbb{R}^{3}} F(u)\,dx -\frac{1}{q}\int_{\mathbb{R}^{3}}|u|^q\,dx,
	\]
	where \(\phi_u^\alpha\) denotes the unique solution of
	\[
	(-\Delta)^\alpha \phi=u^2 \qquad \text{in } \mathbb{R}^3,
	\]
	and
	\[
	F(t)=\int_{0}^{t} s \log s^{2}\,ds =\frac{1}{2} t^{2} \log t^{2}-\frac{t^{2}}{2},
	\]
	is not well defined on \(H^\alpha(\mathbb{R}^3)\), because one may have
	\[
	\int_{\mathbb{R}^{3}}u^{2}\log u^{2}\,dx=-\infty.
	\]
	Following \cite{MR3385171}, we write
	\[
	F_{2}(t)-F_{1}(t)=\frac{1}{2}t^{2}\log t^{2}, \qquad t\in\mathbb{R},
	\]
	where \(F_1\) is a nonnegative convex \(C^1\) function and \(F_2\) is a \(C^1\) function with subcritical growth. This leads naturally to an Orlicz-type framework. More precisely, we consider
	\[
	L^{F_1}(\mathbb{R}^3) = \left\{ u\in L_{\mathrm{loc}}^1(\mathbb{R}^3): \int_{\mathbb{R}^3}F_1(|u|)\,dx<+\infty \right\},
	\]
	and
	\[
	H_{\varepsilon} = \left\{ u\in H^\alpha(\mathbb{R}^3): \int_{\mathbb{R}^3}V(\varepsilon x)|u|^2\,dx<\infty \right\}, \qquad X_{\varepsilon}:=H_{\varepsilon}\cap L^{F_1}(\mathbb{R}^3).
	\]
	In this setting, the energy functional is well defined and belongs to \(C^1\). The use of a Sobolev-Orlicz intersection for scalar logarithmic equations is developed in \cite{AlvesSilva2023}. Here this framework is combined with estimates for the fractional Poisson energy and with a local natural constraint.
	
	A second difficulty is that the potential well is specified on a bounded set, without a compactness assumption at infinity. To recover compactness, we combine the penalization method with an autonomous Nehari-Pohozaev constraint and a local natural constraint near the translated positive ground states. This allows us to recover the topology of
	\[
	M:=\{x\in\Lambda:\ V(x)=V_0\}
	\]
	and then apply Lusternik-Schnirelmann category theory. For \(\delta>0\), write
	\[
	M_\delta:=\{x\in\mathbb R^3:\ \operatorname{dist}(x,M)\le\delta\}.
	\]
	The concentration behavior is obtained from the asymptotic analysis of near-ground-state critical sequences.
	
	We also note that, since
	\[
	|t\log t^2|\le C\bigl(1+|t|^p\bigr), \qquad p\in(2,2_\alpha^*),
	\]
	standard regularity arguments imply that every critical point of the corresponding functional in \(X_\varepsilon\) is a weak solution of the auxiliary problem \eqref{eq2.5}. In particular, every positive critical point satisfying \eqref{eq2.6} yields a positive solution of the original problem \eqref{eq1.3}.
	
	\begin{theorem}\label{Thm1.1}
		Assume that \((V_1)\)-\((V_2)\) hold. Then there exist \(\delta>0\) and \(\varepsilon_{1}>0\) such that, for every \(\varepsilon \in (0,\varepsilon_{1})\), problem \eqref{eq1.2} has at least
		\[
		\operatorname{cat}_{M_{\delta}}(M)
		\]
		distinct positive solutions.
		
		Furthermore, if \((u_\varepsilon,\phi_\varepsilon)\) is one of these distinct positive solutions and \(\eta_{\varepsilon}\in\mathbb{R}^{3}\) is a global maximum point of \(u_\varepsilon\), then
		\[
		\lim_{\varepsilon \to 0} V(\eta_{\varepsilon}) = V_0.
		\]
	\end{theorem}
	Theorem~\ref{Thm1.1} relates the number of positive localized solutions to the topology of the minimum set under a continuous nonnegative potential. Its concentration assertion concerns the solutions constructed in the theorem. The paper is organized as follows. In Section~\ref{Sec2}, we recall some basic facts on Orlicz spaces and fractional Sobolev spaces, and introduce the penalized problem together with the variational setting. In Section~\ref{Sec3}, we prove the existence of a nontrivial solution to the auxiliary problem. In Section~\ref{Sec4}, we study the autonomous Nehari-Pohozaev constraint, construct a local natural constraint for the auxiliary problem, and derive the asymptotic behavior of near-ground-state critical sequences. In Section~\ref{Sec5}, we establish the multiplicity and concentration conclusions in Theorem~\ref{Thm1.1}.
	
	Notation. Throughout the paper,
	\begin{itemize}
		\item \(\|\cdot\|_p\) denotes the norm of \(L^p(\mathbb{R}^3)\), \(p\in[1,\infty]\);
		\item \(o_n(1)\) denotes a quantity converging to \(0\) as \(n\to\infty\);
		\item \(C\) denotes a positive constant whose value may change from line to line.
	\end{itemize}
	
	\section{Preliminaries and Variational Framework}\label{Sec2}
	
	In this section, we recall the material on Orlicz and fractional Sobolev spaces needed later, and introduce the variational setting for problem \eqref{eq1.3}. We refer to \cite{MR2424078,MR2271234,MR1113700,MR2944369,MR3869846,MR4732981,MR3451965,MR3385171} for further background.
	
	\subsection{Orlicz and fractional Sobolev spaces}
	
	\begin{definition}\label{Def2.1}
		A continuous function \(\Phi:\mathbb{R}\to[0,+\infty)\) is called an \(N\)-function if it is convex and even, \(\Phi(t)=0\) if and only if \(t=0\), and
		\[
		\lim_{t\to0^+}\frac{\Phi(t)}{t}=0, \qquad \lim_{t\to+\infty}\frac{\Phi(t)}{t}=+\infty.
		\]
	\end{definition}
	
	We say that \(\Phi\) satisfies the \(\Delta_2\)-condition, and write \(\Phi\in(\Delta_2)\), if there exists \(k>0\) such that
	\[
	\Phi(2t)\le k\,\Phi(t) \qquad \forall t\ge0.
	\]
	Its conjugate function is
	\[
	\widetilde\Phi(s)=\sup_{t\ge0}\{st-\Phi(t)\}, \qquad s\ge0.
	\]
	
	For an open set \(\Omega\subset\mathbb R^N\), the Orlicz space \(L^\Phi(\Omega)\) is endowed with the Luxemburg norm
	\[
	\|u\|_\Phi = \inf\left\{ \lambda>0: \int_\Omega \Phi\!\left(\frac{|u|}{\lambda}\right)\,dx\le1 \right\}.
	\]
	We shall use the standard inequalities
	\[
	st\le \Phi(t)+\widetilde\Phi(s), \qquad s,t\ge0,
	\]
	and
	\[
	\left|\int_\Omega uv\,dx\right| \le 2\|u\|_\Phi\,\|v\|_{\widetilde\Phi}, \qquad u\in L^\Phi(\Omega),\ v\in L^{\widetilde\Phi}(\Omega).
	\]
	If \(\Phi,\widetilde\Phi\in(\Delta_2)\), then \(L^\Phi(\Omega)\) is reflexive and separable, and
	\[
	u_n\to u \ \text{in }L^\Phi(\Omega) \quad\Longleftrightarrow\quad \int_\Omega \Phi(|u_n-u|)\,dx\to0.
	\]
	
	We also recall the standard criterion
	\begin{equation}\label{eq2.1}
		1<l\le \frac{\Phi'(t)t}{\Phi(t)}\le m<\infty, \qquad t\neq0,
	\end{equation}
	which implies \(\Phi,\widetilde\Phi\in(\Delta_2)\). In this case, setting
	\[
	\xi_0(t)=\min\{t^l,t^m\}, \qquad \xi_1(t)=\max\{t^l,t^m\}, \qquad t\ge0,
	\]
	one has
	\begin{equation}\label{eq2.2}
		\xi_0(\|u\|_\Phi) \le \int_\Omega \Phi(|u|)\,dx \le \xi_1(\|u\|_\Phi), \qquad \forall u\in L^\Phi(\Omega).
	\end{equation}
	
	For \(\alpha\in(0,1)\), let \(H^\alpha(\mathbb R^3)=W^{\alpha,2}(\mathbb R^3)\) be the fractional Sobolev space with norm
	\[
	\|u\|_{H^\alpha(\mathbb R^3)}^2 = \int_{\mathbb R^3}\bigl(1+|\xi|^{2\alpha}\bigr)|\mathcal F(u)|^2\,d\xi,
	\]
	where \(\mathcal F\) denotes the Fourier transform. We also denote by \(\mathcal D^{\alpha,2}(\mathbb R^3)\) the completion of \(C_0^\infty(\mathbb R^3)\) with respect to
	\[
	\|u\|_{\mathcal D^{\alpha,2}(\mathbb R^3)} = \left( \iint_{\mathbb R^3\times\mathbb R^3} \frac{|u(x)-u(y)|^2}{|x-y|^{3+2\alpha}}\,dx\,dy \right)^{\frac12} = [u]_{H^\alpha(\mathbb R^3)}.
	\]
	The embedding
	\[
	\mathcal D^{\alpha,2}(\mathbb R^3)\hookrightarrow L^{2_\alpha^*}(\mathbb R^3), \qquad 2_\alpha^*=\frac{6}{3-2\alpha},
	\]
	is continuous, and we denote by
	\[
	S_\alpha := \inf_{u\in\mathcal D^{\alpha,2}(\mathbb R^3)\setminus\{0\}} \frac{\|u\|_{\mathcal D^{\alpha,2}(\mathbb R^3)}^2}{\|u\|_{2_\alpha^*}^2}
	\]
	the best Sobolev constant. The fractional Laplacian is defined in the Fourier sense by
	\[
	\mathcal F\bigl((-\Delta)^\alpha u\bigr)(\xi)=|\xi|^{2\alpha}\mathcal F(u)(\xi),
	\]
	and the usual norms on \(H^\alpha(\mathbb R^3)\) are equivalent.
	
	\subsection{Variational framework}
	
	Set
	\[
	F(t)=\int_0^t s\log s^2\,ds=\frac12 t^2\log t^2-\frac12 t^2.
	\]
	Fix \(\delta\in\bigl(0,e^{-3/2}\bigr)\) sufficiently small that \(M_\delta\subset\Lambda\), and define
	\[
	F_1(s)=
	\begin{cases}
		0, & s=0,\\[0.3em]
		-\frac12 s^2\log s^2, & 0<|s|<\delta,\\[0.3em]
		-\frac12 s^2(\log\delta^2+3)+2\delta |s|-\frac12\delta^2, & |s|\ge \delta,
	\end{cases}
	\]
	and
	\[
	F_2(s)=
	\begin{cases}
		0, & |s|<\delta,\\[0.3em]
		\frac12 s^2\log\!\left(\frac{s^2}{\delta^2}\right)+2\delta |s|-\frac32 s^2-\frac12\delta^2, & |s|\ge \delta.
	\end{cases}
	\]
	Then
	\[
	F_2(s)-F_1(s)=\frac12 s^2\log s^2 \qquad \forall s\in\mathbb R.
	\]
	
	We shall use the following properties:
	\begin{itemize}
		\item[$(f_1)$] \(F_1\) is even, nonnegative, convex, \(F_1\in C^1(\mathbb R,\mathbb R)\), and \(F_1'(s)s\ge0\) for all \(s\in\mathbb R\).
		\item[$(f_2)$] \(F_2\in C^1(\mathbb R,\mathbb R)\cap C^2((\delta,+\infty),\mathbb R)\), and for every \(p\in(2,2_\alpha^*)\) there exists \(C=C_p>0\) such that
		\[
		|F_2'(s)|\le C|s|^{p-1} \qquad \forall s\in\mathbb R.
		\]
		\item[$(f_3)$] The map \(s\mapsto F_2'(s)/s\) is nondecreasing on \((0,+\infty)\) and strictly increasing on \((\delta,+\infty)\).
		\item[$(f_4)$]
		\[
		\lim_{s\to+\infty}\frac{F_2'(s)}{s}=+\infty.
		\]
	\end{itemize}
	
	\begin{proposition}\label{Prop2.1}
		The function \(F_1\) is an \(N\)-function. Moreover, both \(F_1\) and its conjugate \(\widetilde F_1\) satisfy the \(\Delta_2\)-condition.
	\end{proposition}
	
	\begin{proof}
		A direct verification shows that \(F_1\) is even, nonnegative, belongs to \(C^1(\mathbb R)\), vanishes only at \(0\), and
		\[
		\lim_{s\to0}\frac{F_1(s)}{|s|}=0, \qquad \lim_{|s|\to+\infty}\frac{F_1(s)}{|s|}=+\infty.
		\]
		Since \(\delta<e^{-3/2}\), \(F_1\) is convex, hence \(F_1\) is an \(N\)-function. Moreover,
		\[
		\frac{F_1'(s)s}{F_1(s)} = 2+\frac{1}{\log s} \qquad \text{for } 0<s<\delta,
		\]
		and
		\[
		1< \frac{F_1'(s)s}{F_1(s)} \le 2 \qquad \text{for } s\ge\delta.
		\]
		Thus there exists \(l\in(1,2)\) such that
		\begin{equation}\label{eq2.3}
			1<l\le \frac{F_1'(s)s}{F_1(s)}\le 2, \qquad \forall s\neq0.
		\end{equation}
		The conclusion now follows from \eqref{eq2.1}.
	\end{proof}
	
	Choose \(0<a_0<\delta\) so small that
	\[
	b_0:=a_0^{q-2}<\frac{V_0+1}{4},
	\]
	and define
	\[
	\widehat g(t)=
	\begin{cases}
		F_2'(t)+t^{q-1}, & 0\le t\le a_0,\\
		b_0t, & t>a_0.
	\end{cases}
	\]
	Since \(F_2'=0\) on \([0,\delta]\), this function is continuous at \(a_0\). Extend \(\widehat g\) oddly to \(\mathbb R\). Let \(\chi_\Lambda\) be the characteristic function of \(\Lambda\), and set
	\[
	g(x,t)=\chi_\Lambda(x)\bigl(F_2'(t)+|t|^{q-2}t\bigr) +\bigl(1-\chi_\Lambda(x)\bigr)\widehat g(t), \qquad G(x,t)=\int_0^t g(x,s)\,ds.
	\]
	Thus the exterior truncation acts on both \(F_2'\) and the power term. For a fixed \(p\in(2,q)\),
	\begin{equation}\label{eq2.4}
		\begin{aligned}
			\text{(i)}\quad & |g(x,t)|\le b_0|t|+C|t|^{p-1}+|t|^{q-1}, && t\in\mathbb R,\ x\in\mathbb R^3,\\
			\text{(ii)}\quad & 0\le g(x,t)\le F_2'(t)+t^{q-1}, && t\ge0,\ x\in\mathbb R^3,\\
			\text{(iii)}\quad & 0\le g(x,t)t\le b_0t^2, && t\in\mathbb R,\ x\notin\Lambda,\\
			\text{(iv)}\quad & 0\le G(x,t)\le F_2(t)+\frac{|t|^q}{q}, && t\in\mathbb R,\ x\in\mathbb R^3.
		\end{aligned}
	\end{equation}
	Moreover, \(g(x,t)=F_2'(t)+|t|^{q-2}t\) whenever \(|t|\le a_0\), independently of \(x\).
	
	We consider the auxiliary problem
	\begin{equation}\label{eq2.5}
		\begin{cases}
			(-\Delta)^\alpha u+\bigl(V(\varepsilon x)+1\bigr)u+\phi u = g(\varepsilon x,u)-F_1'(u), & \text{in }\mathbb R^3,\\[0.3em]
			(-\Delta)^\alpha \phi=u^2, & \text{in }\mathbb R^3.
		\end{cases}
	\end{equation}
	Let
	\[
	\Lambda_\varepsilon=\{x\in\mathbb R^3:\ \varepsilon x\in\Lambda\}.
	\]
	If \(u\) is a positive solution of \eqref{eq2.5} satisfying
	\begin{equation}\label{eq2.6}
		0<u(x)<a_0 \qquad \forall x\in\mathbb R^3\setminus\Lambda_\varepsilon,
	\end{equation}
	then \((u,\phi_u^\alpha)\) is a positive solution of \eqref{eq1.3}.
	
	For every \(u\in H^\alpha(\mathbb R^3)\), there exists a unique \(\phi_u^\alpha\in\mathcal D^{\alpha,2}(\mathbb R^3)\) solving
	\[
	(-\Delta)^\alpha\phi=u^2 \qquad \text{in }\mathbb R^3,
	\]
	and
	\[
	\phi_u^\alpha(x)=C_\alpha\int_{\mathbb R^3}\frac{u^2(y)}{|x-y|^{3-2\alpha}}\,dy, \qquad x\in\mathbb R^3.
	\]
	We write \(I_{2\alpha}(x)=C_\alpha|x|^{-(3-2\alpha)}\), so that \(\phi_u^\alpha=I_{2\alpha}*u^2\).
	
	Set
	\[
	H_\varepsilon = \left\{ u\in H^\alpha(\mathbb R^3): \int_{\mathbb R^3}V(\varepsilon x)u^2\,dx<+\infty \right\}, \qquad X_\varepsilon=H_\varepsilon\cap L^{F_1}(\mathbb R^3),
	\]
	and equip \(X_\varepsilon\) with the norm
	\[
	\|u\|_\varepsilon=\|u\|_{H_\varepsilon}+\|u\|_{F_1},
	\]
	where
	\[
	\|u\|_{H_\varepsilon} = \left( \int_{\mathbb R^3}\bigl(|(-\Delta)^{\frac{\alpha}{2}}u|^2+(V(\varepsilon x)+1)u^2\bigr)\,dx \right)^{\frac12},
	\]
	and
	\[
	\|u\|_{F_1} = \inf\left\{ \lambda>0: \int_{\mathbb R^3}F_1\!\left(\frac{|u|}{\lambda}\right)\,dx\le1 \right\}.
	\]
	By Proposition~\ref{Prop2.1}, \((X_\varepsilon,\|\cdot\|_\varepsilon)\) is a reflexive and separable Banach space.
	
	Since
	\[
	\alpha\in\left(\frac34,1\right), \qquad 2<\frac{12}{3+2\alpha}<2_\alpha^*,
	\]
	we have
	\[
	H^\alpha(\mathbb R^3)\hookrightarrow L^{\frac{12}{3+2\alpha}}(\mathbb R^3),
	\]
	so the Poisson term is well defined on \(X_\varepsilon\). This embedding already holds for \(\alpha\in[\frac12,1)\), since \(12/(3+2\alpha)\le2_\alpha^*\) precisely when \(\alpha\ge\frac12\). The stronger restriction \(\alpha>\frac34\) makes the interval \(4<q<2_\alpha^*\) nonempty and is also used in the later estimates.
	
	Associated with \eqref{eq2.5}, we define
	\[
	\begin{aligned}
		J_\varepsilon(u) &= \frac12\int_{\mathbb R^3}\bigl(|(-\Delta)^{\frac{\alpha}{2}}u|^2+(V(\varepsilon x)+1)u^2\bigr)\,dx +\frac14\int_{\mathbb R^3}\phi_u^\alpha u^2\,dx\\
		&\quad +\int_{\mathbb R^3}F_1(u)\,dx -\int_{\mathbb R^3}G(\varepsilon x,u)\,dx, \qquad u\in X_\varepsilon.
	\end{aligned}
	\]
	Then \(J_\varepsilon\in C^1(X_\varepsilon,\mathbb R)\), and
	\[
	\begin{aligned}
		\langle J_\varepsilon'(u),v\rangle &= \int_{\mathbb R^3}(-\Delta)^{\frac{\alpha}{2}}u\,(-\Delta)^{\frac{\alpha}{2}}v\,dx +\int_{\mathbb R^3}(V(\varepsilon x)+1)uv\,dx\\
		&\quad +\int_{\mathbb R^3}\phi_u^\alpha uv\,dx +\int_{\mathbb R^3}F_1'(u)v\,dx -\int_{\mathbb R^3}g(\varepsilon x,u)v\,dx,
	\end{aligned}
	\]
	for all \(u,v\in X_\varepsilon\). Hence the critical points of \(J_\varepsilon\) are weak solutions of \eqref{eq2.5}. In particular, every positive critical point satisfying \eqref{eq2.6} yields a positive solution of \eqref{eq1.3}.
	
	We conclude with two standard lemmas.
	
	\begin{lemma}\label{Lem2.1}
		{\rm\cite[Lemma~2.3]{Teng2016}} Let \(u\in H^\alpha(\mathbb R^3)\), with \(\alpha\in\left[\frac12,1\right)\). Then:
		\begin{itemize}
			\item[(i)] \(\phi_u^\alpha\ge0\);
			\item[(ii)] the map \(u\mapsto \phi_u^\alpha\) from \(H^\alpha(\mathbb R^3)\) to \(\mathcal D^{\alpha,2}(\mathbb R^3)\) is continuous and maps bounded sets into bounded sets;
			\item[(iii)]
			\[
			\int_{\mathbb R^3}\phi_u^\alpha u^2\,dx \le C\|u\|_{\frac{12}{3+2\alpha}}^4 \le C\|u\|_{H_\varepsilon}^4 \le C\|u\|_\varepsilon^4;
			\]
			\item[(iv)] \(\phi_{\tau u}^\alpha=\tau^2\phi_u^\alpha\) for all \(\tau\in\mathbb R\);
			\item[(v)] if \(u_n\rightharpoonup u\) in \(H^\alpha(\mathbb R^3)\), then \(\phi_{u_n}^\alpha\rightharpoonup \phi_u^\alpha\) in \(\mathcal D^{\alpha,2}(\mathbb R^3)\);
			\item[(vi)] if \(u_n\to u\) in \(H^\alpha(\mathbb R^3)\), then \(\phi_{u_n}^\alpha\to \phi_u^\alpha\) in \(\mathcal D^{\alpha,2}(\mathbb R^3)\) and
			\[
			\int_{\mathbb R^3}\phi_{u_n}^\alpha u_n^2\,dx \to \int_{\mathbb R^3}\phi_u^\alpha u^2\,dx.
			\]
		\end{itemize}
	\end{lemma}
	
	\begin{lemma}\label{Lem2.2}
		{\rm\cite{MR3216834}} Let \(u\in \mathcal D^{\alpha,2}(\mathbb R^3)\) and \(\varphi\in C_0^\infty(\mathbb R^3)\). For \(r>0\), define
		\[
		\varphi_r(x)=\varphi\!\left(\frac{x}{r}\right).
		\]
		Then:
		\begin{itemize}
			\item[(i)] \(u\varphi_r\to0\) in \(\mathcal D^{\alpha,2}(\mathbb R^3)\) as \(r\to0\);
			\item[(ii)] if \(\varphi\equiv1\) in a neighbourhood of the origin, then
			\[
			u\varphi_r\to u \qquad \text{in }\mathcal D^{\alpha,2}(\mathbb R^3) \qquad \text{as }r\to+\infty.
			\]
		\end{itemize}
	\end{lemma}
	
	\section{Existence of Solutions for the Auxiliary Problem}\label{Sec3}
	
	In this section, we prove the existence of a nontrivial solution to \eqref{eq2.5}. We begin by showing that \(J_\varepsilon\) has the mountain-pass geometry.
	
	\begin{lemma}\label{Lem3.1}
		Let \(\varepsilon>0\) be fixed. Then the functional \(J_\varepsilon\) satisfies:
		\begin{itemize}
			\item[(i)] there exist \(r,\rho>0\) such that
			\[
			J_\varepsilon(u)\ge \rho \qquad \forall\,u\in X_\varepsilon \ \text{with}\ \|u\|_\varepsilon=r;
			\]
			\item[(ii)] there exists \(v\in X_\varepsilon\) with \(\|v\|_\varepsilon>r\) such that
			\[
			J_\varepsilon(v)<0=J_\varepsilon(0).
			\]
		\end{itemize}
	\end{lemma}
	
	\begin{proof}
		Since \(G(x,t)\le F_2(t)+|t|^q/q\) and \(\phi_u^\alpha\ge0\), we have
		\[
		J_\varepsilon(u) \ge \frac12\|u\|_{H_\varepsilon}^2 +\int_{\mathbb R^3}F_1(u)\,dx -\int_{\mathbb R^3}F_2(u)\,dx -\frac1q\int_{\mathbb R^3}|u|^q\,dx.
		\]
		Using \((f_2)\), \eqref{eq2.2}, and the Sobolev embeddings, one obtains
		\[
		J_\varepsilon(u)\ge C_1\|u\|_\varepsilon^2-C_2\|u\|_\varepsilon^p-C_3\|u\|_\varepsilon^q
		\]
		for \(\|u\|_\varepsilon\) small, with \(p>2\) and \(q>4\). This yields (i).
		
		For (ii), choose \(u\in C_0^\infty(\Lambda_\varepsilon)\setminus\{0\}\). Since \(g(\varepsilon x,tu)=F_2'(tu)+|tu|^{q-2}tu\) and \(G(\varepsilon x,tu)=F_2(tu)+|tu|^q/q\) on \(\operatorname{supp}u\), we have
		\[
		J_\varepsilon(tu) = \frac{t^2}{2}\|u\|_{H_\varepsilon}^2 +\frac{t^4}{4}\int_{\mathbb R^3}\phi_u^\alpha u^2\,dx -\frac12\int_{\mathbb R^3}|tu|^2\log |tu|^2\,dx -\frac{t^q}{q}\int_{\mathbb R^3}|u|^q\,dx.
		\]
		Hence
		\[
		J_\varepsilon(tu) = t^2A+t^4B-t^2(\log t)\int_{\mathbb R^3}|u|^2\,dx -\frac{t^q}{q}\int_{\mathbb R^3}|u|^q\,dx
		\]
		for some \(A\in\mathbb R\) and \(B\ge0\). Since \(q>4\), it follows that
		\[
		J_\varepsilon(tu)\to-\infty \qquad \text{as }t\to+\infty.
		\]
		This proves (ii).
	\end{proof}
	
	\begin{lemma}\label{Lem3.2}
		Let \((v_n)\) be a \((PS)_c\)-sequence for \(J_\varepsilon\). Then \((v_n)\) is bounded in \(X_\varepsilon\).
	\end{lemma}
	
	\begin{proof}
		Set
		\[
		A_n=\|v_n\|_{H_\varepsilon}, \qquad B_n=\|v_n\|_{F_1}.
		\]
		Since
		\[
		J_\varepsilon(v_n)\to c, \qquad \|J_\varepsilon'(v_n)\|_{X_\varepsilon'}\to0,
		\]
		we have
		\[
		J_\varepsilon(v_n)-\frac14\langle J_\varepsilon'(v_n),v_n\rangle \le c+1+A_n+B_n
		\]
		for \(n\) large.
		
		We compute the logarithmic contribution to the nonlinear density. From
		\[
		F_2(s)-F_1(s)=\frac12s^2\log s^2, \qquad F_2'(s)s-F_1'(s)s=s^2\log s^2+s^2,
		\]
		we obtain
		\[
		\begin{aligned}
			F_1(s)-\frac14F_1'(s)s-F_2(s)+\frac14F_2'(s)s &=-\frac12s^2\log s^2+\frac14(s^2\log s^2+s^2)\\
			&=\frac14s^2-\frac14s^2\log s^2.
		\end{aligned}
		\]
		The expressions at \(s=0\) are understood by continuity. Consequently,
		\[
		\begin{aligned}
			J_\varepsilon(v_n)-\frac14\langle J_\varepsilon'(v_n),v_n\rangle &= \frac14A_n^2 +\int_{\mathbb R^3} \left( F_1(v_n)-\frac14F_1'(v_n)v_n-G(\varepsilon x,v_n)+\frac14g(\varepsilon x,v_n)v_n \right)\,dx.
		\end{aligned}
		\]
		For \(x\in\Lambda_\varepsilon\), the definition of \(G\) therefore gives
		\[
		\begin{aligned}
			&F_1(s)-\frac14F_1'(s)s-G(\varepsilon x,s)+\frac14g(\varepsilon x,s)s\\
			&\qquad=\frac14s^2-\frac14s^2\log s^2+\frac{q-4}{4q}|s|^q \ge-C_q.
		\end{aligned}
		\]
		Indeed, this expression is continuous at zero and tends to \(+\infty\) as \(|s|\to\infty\), since \(q>4\). On \(\Lambda_\varepsilon^c\), using \eqref{eq2.3} and the definition of \(g\), we obtain
		\[
		F_1(s)-\frac14F_1'(s)s-G(\varepsilon x,s)+\frac14g(\varepsilon x,s)s \ge \frac12F_1(s)-\frac{b_0}{2}|s|^2.
		\]
		Hence, for some \(\kappa>0\),
		\begin{equation}\label{eq3.1}
			\kappa A_n^2+\frac12\int_{\Lambda_\varepsilon^c}F_1(v_n)\,dx \le C_\varepsilon(1+A_n+B_n).
		\end{equation}
		
		Since \(F_1(s)\le C(1+s^2)\),
		\[
		\int_{\Lambda_\varepsilon}F_1(v_n)\,dx\le C_\varepsilon(1+A_n^2),
		\]
		and therefore
		\[
		\int_{\mathbb R^3}F_1(v_n)\,dx \le C_\varepsilon(1+A_n^2) +\int_{\Lambda_\varepsilon^c}F_1(v_n)\,dx.
		\]
		If \(B_n>1\), then \eqref{eq2.2} and \eqref{eq2.3} imply
		\[
		B_n^l\le \int_{\mathbb R^3}F_1(v_n)\,dx
		\]
		for some \(l\in(1,2)\), hence
		\[
		B_n \le C_\varepsilon\Bigl(1+A_n^{2/l} +\Bigl(\int_{\Lambda_\varepsilon^c}F_1(v_n)\,dx\Bigr)^{1/l}\Bigr).
		\]
		Write \(D_n=\int_{\Lambda_\varepsilon^c}F_1(v_n)\,dx\). Since \(2/l<2\) and \(1/l<1\), for every \(\eta>0\),
		\[
		C_\varepsilon\bigl(A_n+A_n^{2/l}+D_n^{1/l}\bigr) \le\eta A_n^2+\eta D_n+C_{\varepsilon,\eta}.
		\]
		Substituting in \eqref{eq3.1} and taking \(\eta\) small enough, we obtain
		\[
		A_n^2+\int_{\Lambda_\varepsilon^c}F_1(v_n)\,dx\le C_\varepsilon.
		\]
		Thus \((A_n)\) is bounded, and so is \(\int_{\mathbb R^3}F_1(v_n)\,dx\). Applying \eqref{eq2.2} again, we infer that \((B_n)\) is bounded. Hence \((v_n)\) is bounded in \(X_\varepsilon\).
	\end{proof}
	
	\begin{lemma}\label{Lem3.3}
		Let \((v_n)\) be a \((PS)_c\)-sequence for \(J_\varepsilon\). Then, for every \(\tau>0\), there exists \(R>0\) such that
		\[
		\limsup_{n\to\infty} \int_{B_R^c(0)} \left( \left|(-\Delta)^{\frac{\alpha}{2}}v_n\right|^2+(V(\varepsilon x)+1)|v_n|^2 \right)\,dx <\tau.
		\]
		Moreover,
		\[
		\lim_{R\to\infty}\limsup_{n\to\infty} \int_{B_R^c(0)}F_1(v_n)\,dx=0.
		\]
	\end{lemma}
	
	\begin{proof}
		Choose \(\eta_R\in C^\infty(\mathbb R^3)\), with \(0\le\eta_R\le1\), equal to zero on \(B_R\), equal to one outside \(B_{2R}\), and satisfying \(|\nabla\eta_R|\le C/R\). Take \(R\) so large that \(\Lambda_\varepsilon\subset B_R\). The sequence \(\eta_R^2v_n\) is bounded in \(X_\varepsilon\), uniformly for \(R\ge1\). We use it as a test function.
		
		For real numbers \(u_x,u_y,\eta_x,\eta_y\),
		\[
		(u_x-u_y)(\eta_x^2u_x-\eta_y^2u_y) =(\eta_xu_x-\eta_yu_y)^2-u_xu_y(\eta_x-\eta_y)^2.
		\]
		Also,
		\[
		\begin{aligned}
			\sup_y\int_{\mathbb R^3} \frac{|\eta_R(x)-\eta_R(y)|^2}{|x-y|^{3+2\alpha}}\,dx &\le\frac{C}{R^2}\int_0^R r^{1-2\alpha}\,dr +C\int_R^\infty r^{-1-2\alpha}\,dr\\
			&\le CR^{-2\alpha}.
		\end{aligned}
		\]
		In particular, the error in the fractional quadratic form satisfies
		\[
		\begin{aligned}
			&\left|\iint_{\mathbb R^3\times\mathbb R^3} \frac{v_n(x)v_n(y)|\eta_R(x)-\eta_R(y)|^2} {|x-y|^{3+2\alpha}}\,dx\,dy\right|\\
			&\quad\le\frac12\iint_{\mathbb R^3\times\mathbb R^3} \frac{\bigl(|v_n(x)|^2+|v_n(y)|^2\bigr)|\eta_R(x)-\eta_R(y)|^2} {|x-y|^{3+2\alpha}}\,dx\,dy\\
			&\quad\le CR^{-2\alpha}\|v_n\|_2^2.
		\end{aligned}
		\]
		Since the Poisson contribution is nonnegative and \(g(\varepsilon x,s)s\le b_0s^2\) on the support of \(\eta_R\), we obtain
		\[
		\begin{aligned}
			&\|(-\Delta)^{\alpha/2}(\eta_Rv_n)\|_2^2 +\int_{\mathbb R^3}(V(\varepsilon x)+1-b_0)\eta_R^2v_n^2\,dx +\int_{\mathbb R^3}\eta_R^2F_1'(v_n)v_n\,dx \le CR^{-2\alpha}+o_n(1).
		\end{aligned}
		\]
		Here \(\varepsilon\) is fixed, and the constant is independent of \(n,R\). By \eqref{eq2.3}, this gives both the weighted \(L^2\) tail estimate and the stated Orlicz tail estimate.
		
		To obtain the tail of \((-\Delta)^{\alpha/2}v_n\), note that the kernel of the commutator \([(-\Delta)^{\alpha/2},\eta_R]\) is bounded in absolute value by
		\[
		C\frac{\min\{1,|x-y|/R\}}{|x-y|^{3+\alpha}}.
		\]
		Its integrals in either variable are bounded by \(CR^{-\alpha}\), since \(\alpha<1\). Schur's estimate therefore gives
		\[
		\|[(-\Delta)^{\alpha/2},\eta_R]v_n\|_2 \le CR^{-\alpha}\|v_n\|_2.
		\]
		Combining this with the preceding estimate for \(\eta_Rv_n\) proves the first assertion.
	\end{proof}
	
	\begin{lemma}\label{Lem3.4}
		The functional \(J_\varepsilon\) satisfies the \((PS)\)-condition.
	\end{lemma}
	
	\begin{proof}
		Let \((v_n)\) be a \((PS)_c\)-sequence for \(J_\varepsilon\). By Lemma~\ref{Lem3.2}, \((v_n)\) is bounded in \(X_\varepsilon\). Hence, up to a subsequence,
		\[
		v_n\rightharpoonup v \quad \text{in }X_\varepsilon, \qquad v_n\rightharpoonup v \quad \text{in }H^\alpha(\mathbb R^3), \qquad v_n\to v \quad \text{a.e. in }\mathbb R^3.
		\]
		By the local compact embedding and Lemma~\ref{Lem3.3}, we have
		\[
		v_n\to v \qquad \text{in }L^r(\mathbb R^3) \quad \forall r\in[2,2_\alpha^*).
		\]
		In particular,
		\begin{equation}\label{eq3.2}
			v_n\to v \qquad \text{in }L^2(\mathbb R^3)\cap L^p(\mathbb R^3)\cap L^q(\mathbb R^3).
		\end{equation}
		
		For compactly supported smooth test functions, local compactness and the growth bounds for \(F_1'\) allow passage to the limit in the logarithmic term. Lemma~\ref{Lem2.1} and \eqref{eq3.2} handle the remaining nonlinear terms. Density then gives
		\[
		J_\varepsilon'(v)=0.
		\]
		Set
		\[
		d_n:=v_n-v.
		\]
		Then
		\[
		d_n\rightharpoonup0 \quad \text{in }X_\varepsilon, \qquad d_n\to0 \quad \text{in }L^r(\mathbb R^3) \quad \forall r\in[2,2_\alpha^*).
		\]
		
		By Lemma~\ref{Lem2.1} and Hölder's inequality,
		\[
		\int_{\mathbb R^3}\bigl(\phi_{v_n}^\alpha v_n-\phi_v^\alpha v\bigr)d_n\,dx\to0.
		\]
		Moreover, by \eqref{eq2.4}(i) and \eqref{eq3.2},
		\[
		\int_{\mathbb R^3}\bigl(g(\varepsilon x,v_n)-g(\varepsilon x,v)\bigr)d_n\,dx\to0.
		\]
		Hence
		\[
		\langle J_\varepsilon'(v_n)-J_\varepsilon'(v),d_n\rangle=o_n(1)
		\]
		gives
		\[
		\|d_n\|_{H_\varepsilon}^2 + \int_{\mathbb R^3}\bigl(F_1'(v_n)-F_1'(v)\bigr)d_n\,dx \to0.
		\]
		Since \(F_1\) is convex, the second term is nonnegative, and therefore
		\[
		v_n\to v \qquad \text{in }H_\varepsilon.
		\]
		
		It remains to prove convergence of the Orlicz modular. On each fixed ball, for sufficiently small \(\sigma>0\),
		\[
		F_1(s)\le C_\sigma\bigl(|s|^{2-\sigma}+|s|^{2+\sigma}\bigr),
		\]
		so local compactness gives \(\int_{B_R}F_1(d_n)\,dx\to0\). Outside the ball, the \(\Delta_2\)-condition yields
		\[
		F_1(d_n)\le C\bigl(F_1(v_n)+F_1(v)\bigr).
		\]
		Lemma~\ref{Lem3.3} and Fatou's lemma make these tails uniformly small. Letting first \(n\to\infty\) and then \(R\to\infty\), we obtain
		\[
		\int_{\mathbb R^3}F_1(d_n)\,dx\to0.
		\]
		By \eqref{eq2.2} and \eqref{eq2.3}, \(v_n\to v\) in \(L^{F_1}(\mathbb R^3)\). Thus
		\[
		v_n\to v \qquad \text{in }X_\varepsilon.
		\]
	\end{proof}
	
	\begin{theorem}\label{Thm3.1}
		For each \(\varepsilon>0\), the functional \(J_\varepsilon\) has a nontrivial nonnegative critical point \(u_\varepsilon\in X_\varepsilon\). Consequently, \eqref{eq2.5} admits a nontrivial nonnegative weak solution.
	\end{theorem}
	
	\begin{proof}
		For the construction of a nonnegative solution, define
		\[
		J_\varepsilon^+(u) = \frac12\|u\|_{H_\varepsilon}^2 +\frac14\int_{\mathbb R^3}\phi_u^\alpha u^2\,dx +\int_{\mathbb R^3}F_1(u)\,dx -\int_{\mathbb R^3}G(\varepsilon x,u^+)\,dx.
		\]
		The function \(s\mapsto G(x,s^+)\) has continuous derivative equal to \(g(x,s)\) for \(s>0\) and zero for \(s\le0\). Thus \(J_\varepsilon^+\in C^1(X_\varepsilon,\mathbb R)\). The full convex term \(F_1(u)\) is retained.
		
		By the proof of Lemma~\ref{Lem3.1}, \(J_\varepsilon^+\) has the mountain-pass geometry: the lower bound near zero is unchanged, and a nonnegative function in \(C_c^\infty(\Lambda_\varepsilon)\setminus\{0\}\) gives a path of negative energy. Set
		\[
		\mathscr C_\varepsilon^+ = \left\{ \gamma\in C([0,1],X_\varepsilon): \gamma(0)=0,\ J_\varepsilon^+(\gamma(1))<0 \right\},
		\]
		and
		\[
		\mathfrak c_\varepsilon^+ = \inf_{\gamma\in\mathscr C_\varepsilon^+} \max_{t\in[0,1]}J_\varepsilon^+(\gamma(t))>0.
		\]
		The estimates of Lemmas~\ref{Lem3.2}-\ref{Lem3.4} also hold for \(J_\varepsilon^+\). Indeed, on the positive set they are identical; on the negative set the nonlinear contribution to \(J_\varepsilon^+(u)-\frac14(J_\varepsilon^+)'(u)u\) is
		\[
		F_1(u)-\frac14F_1'(u)u\ge\frac12F_1(u).
		\]
		The same cutoff test gives the exterior estimates, and the convexity argument gives strong convergence. Consequently, the mountain-pass theorem gives \(u_\varepsilon\in X_\varepsilon\) such that
		\[
		(J_\varepsilon^+)'(u_\varepsilon)=0, \qquad J_\varepsilon^+(u_\varepsilon)=\mathfrak c_\varepsilon^+>0.
		\]
		In particular, \(u_\varepsilon\ne0\).
		
		Write \(u_\varepsilon^-=\max\{-u_\varepsilon,0\}\). Testing with \(-u_\varepsilon^-\), the fractional quadratic form satisfies
		\[
		\int_{\mathbb R^3}(-\Delta)^{\alpha/2}u_\varepsilon (-\Delta)^{\alpha/2}(-u_\varepsilon^-)\,dx \ge \|(-\Delta)^{\alpha/2}u_\varepsilon^-\|_2^2.
		\]
		The truncated positive nonlinearity vanishes on this test. Since \(\phi_{u_\varepsilon}^\alpha\ge0\), it follows that
		\[
		\|u_\varepsilon^-\|_{H_\varepsilon}^2 +\int_{\mathbb R^3}F_1'(u_\varepsilon^-)u_\varepsilon^-\,dx \le0.
		\]
		Thus \(u_\varepsilon\ge0\). At a nonnegative function the derivatives of \(J_\varepsilon^+\) and \(J_\varepsilon\) agree. Therefore
		\[
		J_\varepsilon'(u_\varepsilon)=0, \qquad J_\varepsilon(u_\varepsilon)=\mathfrak c_\varepsilon^+>0,
		\]
		and \(u_\varepsilon\) is the required nontrivial nonnegative weak solution.
	\end{proof}
	
	\begin{lemma}\label{Lem3.5}
		Every nontrivial nonnegative critical point of \(J_\varepsilon\) is strictly positive in \(\mathbb R^3\).
	\end{lemma}
	
	\begin{proof}
		Let \(u\in X_\varepsilon\) be a nontrivial nonnegative critical point of \(J_\varepsilon\). Since \(V\ge0\), \(\phi_u^\alpha\ge0\), and \(F_1'(s)\ge0\) for \(s\ge0\), \eqref{eq2.4} gives
		\[
		\bigl(( -\Delta)^\alpha+1\bigr)u \le C u+C u^{q-1}
		\]
		in the distributional sense. Here and below the constants in this estimate are independent of \(\varepsilon\).
		
		Let \(K_\alpha\) be the positive convolution kernel with Fourier transform \((1+|\xi|^{2\alpha})^{-1}\). By \cite[Theorem~3.3]{FQT},
		\[
		K_\alpha\in L^a(\mathbb R^3) \qquad \text{for }1\le a<\frac3{3-2\alpha}.
		\]
		In particular, \(K_\alpha\in L^1\cap L^2\), since \(4\alpha>3\). The weak comparison principle for \(( -\Delta)^\alpha+1\) yields
		\[
		0\le u\le C K_\alpha*u+C K_\alpha*u^{q-1}.
		\]
		This comparison follows by applying the positive resolvent to the nonnegative right-hand side and testing the difference with its positive part; truncation and density justify the test.
		
		We give the integrability iteration, since its uniformity will be used later. Put \(p_0=2_\alpha^*\), so that \(u\in L^2\cap L^{p_0}\). Choose \(\sigma>0\) such that
		\[
		\frac{q-2}{p_0}<\frac{2\alpha}{3}-\sigma<\frac{2\alpha}{3}.
		\]
		Such a choice is possible precisely because \(q<2_\alpha^*\). The term \(K_\alpha*u\) is bounded in \(L^2\cap L^\infty\), and hence in every \(L^p\), \(2\le p\le\infty\). For the other term, Young's inequality gives, as long as the right-hand side is positive,
		\[
		\frac1{p_{j+1}} = \frac{q-1}{p_j}-\left(\frac{2\alpha}{3}-\sigma\right), \qquad \|K_\alpha*u^{q-1}\|_{p_{j+1}} \le C\|u\|_{p_j}^{q-1}.
		\]
		The reciprocal exponents strictly decrease. The recurrence reaches, after finitely many steps, an exponent \(p_j\) for which
		\[
		\frac{p_j}{q-1}>\frac3{2\alpha}.
		\]
		At this stage H\"older's inequality in the convolution yields an \(L^\infty\) bound. Intermediate exponents may be chosen to avoid an equality at an endpoint. It follows that
		\[
		u\in L^\infty(\mathbb R^3),
		\]
		and the bound is uniform for every family of nonnegative solutions bounded in \(H^\alpha(\mathbb R^3)\).
		
		Splitting the Riesz integral into \(|x-y|<1\) and \(|x-y|\ge1\), we obtain
		\[
		\|\phi_u^\alpha\|_\infty \le C\bigl(\|u\|_\infty^2+\|u\|_2^2\bigr).
		\]
		Since \(V\) is bounded on compact sets and \(s\log s^2\) is bounded on bounded intervals, the right-hand side of the equation is locally bounded. Interior regularity therefore gives
		\[
		u\in C_{\mathrm{loc}}^{0,\sigma_0}(\mathbb R^3)
		\]
		for some \(\sigma_0\in(0,1)\).
		
		Suppose that \(u(x_0)=0\). By continuity, \(u<a_0\) on a small ball about \(x_0\), so the original and the penalized nonlinearities agree there. The continuous weak solution satisfies the viscosity equation on this ball by \cite{ServadeiValdinoci}. Use the zero function as an interior touching test and retain \(u\) outside a smaller ball. Its fractional Laplacian at \(x_0\) is
		\[
		-C_\alpha\int_{\mathbb R^3\setminus B_r(x_0)} \frac{u(y)}{|x_0-y|^{3+2\alpha}}\,dy<0
		\]
		for \(r\) sufficiently small, because \(u\not\equiv0\). All local terms in the equation vanish at \(x_0\), which contradicts the supersolution inequality. Thus \(u>0\) throughout \(\mathbb R^3\).
	\end{proof}
	
	\section{The Nehari-Pohozaev Constraint and Positive Solutions}\label{Sec4}
	
	In this section, we study the autonomous problem on a Nehari-Pohozaev constraint, establish the existence of a positive ground state, and construct a local natural constraint for the auxiliary problem. We then analyze the asymptotic behavior of its low-energy solutions.
	
	We first establish a lower bound for nontrivial critical points. Set
	\[
	\mathcal N_{\varepsilon} := \left\{ u\in X_\varepsilon\setminus\{0\}:J_\varepsilon'(u)u=0 \right\}.
	\]
	Every nontrivial critical point of \(J_\varepsilon\) belongs to this set.
	
	\begin{lemma}\label{Lem4.1}
		There exists \(\beta_0>0\) such that, for all \(\varepsilon>0\),
		\[
		\|u\|_{\varepsilon}\ge \|u\|_{H_{\varepsilon}}\ge \beta_0, \qquad \forall\,u\in \mathcal N_{\varepsilon}.
		\]
	\end{lemma}
	
	\begin{proof}
		Let \(u\in\mathcal N_\varepsilon\). Since \(J_\varepsilon'(u)u=0\), by \eqref{eq2.4}(i), \(\phi_u^\alpha\ge0\), and \(F_1'(s)s\ge0\), we obtain
		\begin{equation}\label{eq4.1}
			\begin{aligned}
				0 &= \|u\|_{H_\varepsilon}^{2} +\int_{\mathbb R^3}\phi_u^{\alpha}u^2\,dx +\int_{\mathbb R^3}F_1'(u)u\,dx -\int_{\mathbb R^3}g(\varepsilon x,u)u\,dx\\
				&\ge \|u\|_{H_\varepsilon}^{2} -b_0\int_{\mathbb R^3}|u|^2\,dx -C\int_{\mathbb R^3}|u|^p\,dx -\int_{\mathbb R^3}|u|^q\,dx\\
				&\ge \frac12\,\|u\|_{H_\varepsilon}^{2} -C_1\|u\|_{H_\varepsilon}^{p} -C_2\|u\|_{H_\varepsilon}^{q},
			\end{aligned}
		\end{equation}
		where we used \(b_0<(V_0+1)/4\) and the Sobolev embeddings. Since \(p>2\) and \(q>4\), \eqref{eq4.1} yields the existence of \(\beta_0>0\), independent of \(\varepsilon\), such that
		\[
		\|u\|_{H_\varepsilon}\ge \beta_0.
		\]
		The inequality \(\|u\|_{\varepsilon}\ge \|u\|_{H_\varepsilon}\) follows directly from the definition of \(\|\cdot\|_\varepsilon\).
	\end{proof}
	
	\subsection{The autonomous problem}
	
	We now consider the autonomous problem
	\begin{equation}\label{eq4.2}
		\begin{cases}
			(-\Delta)^\alpha u+V_0\,u+\phi\,u=u \log u^2+|u|^{q-2}u, & \text{in }\mathbb{R}^3,\\[2pt]
			(-\Delta)^\alpha \phi=u^2, & \text{in }\mathbb{R}^3.
		\end{cases}
	\end{equation}
	Its associated functional is
	\[
	\begin{aligned}
		J_{0}(u) &= \frac{1}{2}\int_{\mathbb{R}^{3}}\bigl(|(-\Delta)^{\frac{\alpha}{2}} u|^{2}+(V_{0}+1)|u|^{2}\bigr)\,dx +\frac{1}{4}\int_{\mathbb{R}^{3}}\phi_u^{\alpha}\,u^{2}\,dx\\
		&\quad +\int_{\mathbb{R}^{3}} F_{1}(u)\,dx -\int_{\mathbb{R}^{3}} F_{2}(u)\,dx -\frac{1}{q}\int_{\mathbb{R}^{3}}|u|^q\,dx .
	\end{aligned}
	\]
	Let
	\[
	X:=H^\alpha(\mathbb{R}^{3})\cap L^{F_{1}}(\mathbb{R}^{3}), \qquad \|u\|_X:=\|u\|_{H_0}+\|u\|_{F_1},
	\]
	and
	\[
	\|u\|_{H_0}^2:=\int_{\mathbb{R}^{3}}\bigl(|(-\Delta)^{\frac{\alpha}{2}}u|^2+(V_0+1)u^2\bigr)\,dx.
	\]
	
	By \eqref{eq2.4}, for every \(u\in X_\varepsilon\),
	\begin{equation}\label{eq4.3}
		\begin{aligned}
			J_\varepsilon(u)-J_0(u) &=\frac12\int_{\mathbb R^3}(V(\varepsilon x)-V_0)u^2\,dx\\
			&\quad+\int_{\mathbb R^3} \left(F_2(u)+\frac{|u|^q}{q}-G(\varepsilon x,u)\right)\,dx \ge0.
		\end{aligned}
	\end{equation}
	Both integrands on the right are nonnegative. This comparison will also be used for functions on the local constraint below.
	
	We record two estimates used in the autonomous compactness argument.
	
	\begin{lemma}\label{Lem4.2}
		For every \(\kappa\in(0,1)\) and every fixed \(r\in(0,1)\), there exists \(C_{\kappa,r}>0\) such that
		\[
		\int_{\mathbb R^3}u^2\log u^2\,dx \le \kappa\,\|u\|_{H_\varepsilon}^2 + C_{\kappa,r}\bigl(1+\|u\|_{H_\varepsilon}^{2(1+r)}\bigr) \qquad \forall u\in H_\varepsilon.
		\]
		The constant \(C_{\kappa,r}\) is independent of \(\varepsilon\).
	\end{lemma}
	
	\begin{proof}
		By the fractional logarithmic Sobolev inequality and the estimate
		\[
		\|u\|_2^2\le \frac1{V_0+1}\|u\|_{H_\varepsilon}^2,
		\]
		one obtains
		\[
		\int_{\mathbb R^3}u^2\log u^2\,dx \le \kappa\,\|u\|_{H_\varepsilon}^2 + C_1\|u\|_2^2\bigl|\log \|u\|_2^2\bigr| + C_2.
		\]
		The conclusion follows from the elementary inequality \(t|\log t|\le C_r(1+t^{1+r})\), for \(t>0\).
	\end{proof}
	
	\begin{lemma}\label{Lem4.3}
		Let \((u_n)\subset X\) be bounded and suppose that for some \(R>0\),
		\[
		\sup_{y\in\mathbb R^3}\int_{B_R(y)}|u_n|^2\,dx \to0.
		\]
		Then
		\[
		u_n\to0 \qquad \text{in }L^r(\mathbb R^3), \quad \forall r\in(2,2_\alpha^*).
		\]
	\end{lemma}
	
	\begin{proof}
		This is the standard Lions vanishing lemma in the fractional Sobolev setting.
	\end{proof}
	
	For the following scaling computations, write
	\[
	\begin{gathered}
		A(u)=\|(-\Delta)^{\alpha/2}u\|_2^2, \qquad m(u)=\|u\|_2^2, \qquad B(u)=\int_{\mathbb R^3}\phi_u^\alpha u^2\,dx,\\
		L(u)=\int_{\mathbb R^3}u^2\log u^2\,dx, \qquad Q(u)=\|u\|_q^q.
	\end{gathered}
	\]
	The integral defining \(L(u)\) is finite for \(u\in X\), by the splitting \(s^2\log s^2=2(F_2-F_1)(s)\). For \(t>0\), put
	\[
	T_tu(x)=t^3u(tx).
	\]
	Changing variables in each term gives
	\[
	\begin{aligned}
		A(T_tu)&=t^{3+2\alpha}A(u), & m(T_tu)&=t^3m(u),\\
		B(T_tu)&=t^{9-2\alpha}B(u), & Q(T_tu)&=t^{3(q-1)}Q(u),\\
		L(T_tu)&=t^3L(u)+6t^3\log t\,m(u).
	\end{aligned}
	\]
	For the Poisson term, the equality follows directly from
	\[
	B(T_tu) =C_\alpha t^{12} \iint_{\mathbb R^3\times\mathbb R^3} \frac{u(tx)^2u(ty)^2}{|x-y|^{3-2\alpha}}\,dx\,dy.
	\]
	Thus
	\begin{equation}\label{eq4.4}
		\begin{aligned}
			E_u(t):=J_0(T_tu) &=\frac{t^{3+2\alpha}}2A(u) +\frac{t^{9-2\alpha}}4B(u) +\frac{t^3}2\bigl((V_0+1)m(u)-L(u)\bigr)\\
			&\quad-3t^3\log t\,m(u)-\frac{t^{3(q-1)}}qQ(u).
		\end{aligned}
	\end{equation}
	The scaling is strongly continuous on \((0,\infty)\times X\). Indeed, the Sobolev part follows by density, while the Orlicz part follows by approximation with bounded compactly supported functions and the dilation estimates implied by \eqref{eq2.3}. For \(0<t\le1\), these estimates give
	\[
	\int_{\mathbb R^3}F_1(T_tu)\,dx =t^{-3}\int_{\mathbb R^3}F_1(t^3u)\,dx \le t^{3(l-1)}\int_{\mathbb R^3}F_1(u)\,dx.
	\]
	Hence \(T_tu\to0\) in \(X\) as \(t\downarrow0\).
	
	Define
	\[
	\mathcal G_0(u) =\left.\frac d{dt}J_0(T_tu)\right|_{t=1}, \qquad \mathcal P_0=\{u\in X\setminus\{0\}:\mathcal G_0(u)=0\}, \qquad c_0=\inf_{u\in\mathcal P_0}J_0(u).
	\]
	Differentiation of \eqref{eq4.4} yields
	\begin{equation}\label{eq4.5}
		\begin{aligned}
			\mathcal G_0(u) &=\frac{3+2\alpha}{2}A(u) +\frac{9-2\alpha}{4}B(u) +\frac32(V_0-1)m(u)\\
			&\quad-\frac32L(u)-\frac{3(q-1)}qQ(u).
		\end{aligned}
	\end{equation}
	In particular, \(\mathcal G_0\in C^1(X,\mathbb R)\), by the Orlicz and Sobolev estimates in Section~\ref{Sec2}.
	
	\begin{lemma}\label{Lem4.4}
		For every \(u\in X\setminus\{0\}\), there exists a unique \(t(u)>0\) such that
		\[
		T_{t(u)}u\in\mathcal P_0.
		\]
		Moreover,
		\[
		J_0(T_{t(u)}u)=\max_{t>0}J_0(T_tu),
		\]
		the maximum is strict, and \(u\mapsto t(u)\) is continuous on \(X\setminus\{0\}\).
	\end{lemma}
	
	\begin{proof}
		Put
		\[
		H_u(t)=\frac{E_u'(t)}{t^{8-2\alpha}}.
		\]
		A direct differentiation of \eqref{eq4.4} gives
		\begin{equation}\label{eq4.6}
			\begin{aligned}
				H_u'(t) &=-\frac{(3+2\alpha)(6-4\alpha)}2 A(u)t^{4\alpha-7}\\
				&\quad+3(3-\alpha)t^{2\alpha-7} \int_{\mathbb R^3}u^2 \left[ \log(t^6u^2)-V_0-\frac{\alpha}{3-\alpha} -C_{\alpha,q}(t^3|u|)^{q-2} \right]\,dx,
			\end{aligned}
		\end{equation}
		where
		\[
		C_{\alpha,q} =\frac{(q-1)\bigl(3(q-4)+2\alpha\bigr)}{q(3-\alpha)}.
		\]
		Both \((q-1)/q\) and \(3(q-4)+2\alpha\) increase for \(q>4\). Consequently,
		\[
		C_{\alpha,q} >\frac{3\alpha}{2(3-\alpha)}>\frac12,
		\]
		because \(\alpha>3/4\). For every \(s>0\), maximizing \(\log s^2/s^{q-2}\) gives
		\[
		\log s^2\le\frac{2}{e(q-2)}s^{q-2} <\frac12s^{q-2}.
		\]
		Since \(V_0\ge0\), the integrand in brackets in \eqref{eq4.6} is at most \(-\alpha/(3-\alpha)\). Therefore
		\[
		H_u'(t)<0 \qquad \forall t>0.
		\]
		
		As \(t\downarrow0\), the positive contribution \(-9m(u)t^{2\alpha-6}\log t\) to \(H_u(t)\) dominates the other terms of the same or smaller order, so
		\[
		\lim_{t\downarrow0}H_u(t)=+\infty.
		\]
		At infinity, the power contribution is
		\[
		-\frac{3(q-1)}qQ(u)t^{3(q-4)+2\alpha},
		\]
		where \(Q(u)>0\). The remaining nonconstant terms tend to zero, and the Poisson contribution is constant. Hence
		\[
		\lim_{t\to\infty}H_u(t)=-\infty.
		\]
		There is therefore exactly one zero \(t(u)\). Since
		\[
		\mathcal G_0(T_tu)=tE_u'(t),
		\]
		this is the unique point on the scaling curve belonging to \(\mathcal P_0\). The sign change of \(H_u\) proves that it is the unique global maximum of \(E_u\).
		
		Finally, let \(u_n\to u\ne0\) in \(X\). Each coefficient in \eqref{eq4.4} is continuous in \(u\). Thus \(H_{u_n}(a)\to H_u(a)\) for each \(a>0\). Given \(0<a<t(u)<b\), we have \(H_u(a)>0>H_u(b)\), and the same inequalities hold for \(u_n\) when \(n\) is large. Hence \(a<t(u_n)<b\). Letting \(a\uparrow t(u)\) and \(b\downarrow t(u)\) proves \(t(u_n)\to t(u)\).
	\end{proof}
	
	Set
	\[
	\mathcal K_0(u)=J_0(u)-\frac{\mathcal G_0(u)}{9-2\alpha}.
	\]
	Subtracting \eqref{eq4.5} from the energy, we obtain
	\begin{equation}\label{eq4.7}
		\begin{aligned}
			\mathcal K_0(u) &=\frac{3-2\alpha}{9-2\alpha}A(u) +\frac{(3-\alpha)V_0+6-\alpha}{9-2\alpha}m(u)\\
			&\quad-\frac{3-\alpha}{9-2\alpha}L(u) +\frac{3(q-4)+2\alpha}{q(9-2\alpha)}Q(u).
		\end{aligned}
	\end{equation}
	
	\begin{lemma}\label{Lem4.5}
		There exist \(\kappa_0,\rho_0>0\) such that
		\begin{equation}\label{eq4.8}
			\mathcal K_0(u)\ge \kappa_0\bigl(\|u\|_{H^\alpha}^2+\|u\|_q^q\bigr) \qquad \forall u\in X,
		\end{equation}
		and
		\begin{equation}\label{eq4.9}
			\mathcal G_0(u)>0 \quad\text{if }0<\|u\|_{H^\alpha}\le\rho_0.
		\end{equation}
		For every \(u\ne0\), the function \(t\mapsto\mathcal K_0(T_tu)\) is strictly increasing. In particular, \(c_0>0\), and every minimizing sequence for \(c_0\) on \(\mathcal P_0\) is bounded in \(X\).
	\end{lemma}
	
	\begin{proof}
		For all \(s>0\),
		\[
		\log s^2\le1+\frac16s^{q-2}.
		\]
		For \(s\le1\) this is immediate. For \(s\ge1\), put \(y=s^{q-2}\) and use \(2/(q-2)<1\) together with
		\[
		\log y-\frac y6\le\log6-1<1.
		\]
		Substituting in \eqref{eq4.7} gives
		\[
		\begin{aligned}
			\mathcal K_0(u) &\ge\frac{3-2\alpha}{9-2\alpha}A(u) +\frac{(3-\alpha)V_0+3}{9-2\alpha}m(u)\\
			&\quad+\frac1{9-2\alpha} \left(\frac{3(q-4)+2\alpha}{q}-\frac{3-\alpha}{6}\right)Q(u).
		\end{aligned}
		\]
		The last coefficient is positive, since
		\[
		\frac{3(q-4)+2\alpha}{q}-\frac{3-\alpha}{6} >\frac{4\alpha-3}{6}>0.
		\]
		The other two coefficients are also positive. This proves \eqref{eq4.8}.
		
		To prove \eqref{eq4.9}, fix \(p\in(2,q)\). Since \(-\log s^2\to+\infty\) as \(s\to0\), while the logarithm grows more slowly than \(s^{p-2}\) at infinity, there is \(C_p>0\) such that
		\[
		\frac32(V_0-1)s^2-\frac32s^2\log s^2 \ge s^2-C_p|s|^p.
		\]
		Using this in \eqref{eq4.5}, dropping the nonnegative Poisson term, and applying the Sobolev embeddings, we find
		\[
		\mathcal G_0(u) \ge c\|u\|_{H^\alpha}^2 -C_1\|u\|_{H^\alpha}^p-C_2\|u\|_{H^\alpha}^q.
		\]
		This is positive for sufficiently small nonzero \(\|u\|_{H^\alpha}\).
		
		Furthermore, since \(\mathcal G_0(T_tu)=tE_u'(t)\),
		\begin{equation}\label{eq4.10}
			\frac d{dt}\mathcal K_0(T_tu) =-\frac{t^{9-2\alpha}}{9-2\alpha}H_u'(t)>0.
		\end{equation}
		The set \(\mathcal P_0\) is nonempty by Lemma~\ref{Lem4.4}. On it, \(J_0=\mathcal K_0\) and \(\|u\|_{H^\alpha}>\rho_0\). Hence
		\[
		c_0\ge\kappa_0\rho_0^2>0.
		\]
		
		Let \(u_n\in\mathcal P_0\) satisfy \(J_0(u_n)\to c_0\). Then \eqref{eq4.8} bounds \(u_n\) in \(H^\alpha\) and in \(L^q\). Lemma~\ref{Lem4.2} also bounds \(L(u_n)\) from above. To control the negative logarithmic part, substitute \(L=2\int(F_2-F_1)\) into \eqref{eq4.7}. All terms except \(-2(3-\alpha)\int F_2(u_n)/(9-2\alpha)\) are nonnegative, so
		\[
		\frac{2(3-\alpha)}{9-2\alpha} \int_{\mathbb R^3}F_1(u_n)\,dx \le c_0+1+C\int_{\mathbb R^3}F_2(u_n)\,dx\le C.
		\]
		The last inequality follows from \((f_2)\) and the Sobolev embedding. Equations \eqref{eq2.2}-\eqref{eq2.3} now bound \(\|u_n\|_{F_1}\). Thus \((u_n)\) is bounded in \(X\).
	\end{proof}
	
	\begin{proposition}\label{Prop4.1}
		Every minimizing sequence \((u_n)\subset\mathcal P_0\) for \(c_0\) has a subsequence and translations \(y_n\in\mathbb R^3\) such that
		\[
		u_n(\cdot+y_n)\to u_0 \qquad\text{in }X,
		\]
		where
		\[
		u_0\in\mathcal P_0, \qquad J_0(u_0)=c_0, \qquad J_0'(u_0)=0.
		\]
		Every minimizer has a fixed sign. Moreover,
		\begin{equation}\label{eq4.11}
			c_0 =\inf_{u\in X\setminus\{0\}}\max_{t>0}J_0(T_tu) =\inf_{\gamma\in\mathscr C_0}\max_{s\in[0,1]}J_0(\gamma(s)),
		\end{equation}
		where
		\[
		\mathscr C_0 =\{\gamma\in C([0,1],X):\gamma(0)=0,\ J_0(\gamma(1))<0\}.
		\]
	\end{proposition}
	
	\begin{proof}
		We first establish \eqref{eq4.11}. The first equality follows from Lemma~\ref{Lem4.4}. By \eqref{eq4.8}, negative energy implies \(\mathcal G_0<0\). On the other hand, \eqref{eq4.9} shows that \(\mathcal G_0\) is positive when the Sobolev norm is small and nonzero. Each path in \(\mathscr C_0\) therefore meets \(\mathcal P_0\), and its maximum energy is at least \(c_0\). Conversely, \(T_tu\to0\) in \(X\) as \(t\downarrow0\), and \eqref{eq4.4} gives \(J_0(T_tu)\to-\infty\) as \(t\to\infty\). Stopping this scaling curve at a sufficiently large parameter gives a path in \(\mathscr C_0\). This proves the reverse inequality.
		
		Let \(u_n\in\mathcal P_0\) and \(J_0(u_n)\to c_0\). Lemma~\ref{Lem4.5} bounds this sequence in \(X\). It cannot vanish in the sense of Lemma~\ref{Lem4.3}. Indeed, vanishing would give
		\[
		u_n\to0\quad\text{in }L^p(\mathbb R^3)\cap L^q(\mathbb R^3)
		\]
		for a fixed \(p\in(2,q)\). The estimate in the proof of \eqref{eq4.9}, together with \(\mathcal G_0(u_n)=0\), would then imply \(\|u_n\|_{H^\alpha}\to0\), contrary to \eqref{eq4.9}. Thus there exist \(R>0\), \(\gamma>0\), and \(y_n\in\mathbb R^3\) such that
		\[
		\int_{B_R(y_n)}|u_n|^2\,dx\ge\gamma.
		\]
		Define
		\[
		v_n(x)=u_n(x+y_n).
		\]
		By translation invariance, \(v_n\in\mathcal P_0\) and \(J_0(v_n)\to c_0\). Passing to a subsequence, we have
		\[
		\begin{gathered}
			v_n\rightharpoonup u_0\quad\text{in }X, \qquad v_n\to u_0\quad\text{a.e. in }\mathbb R^3,\\
			v_n\to u_0\quad\text{in }L^p_{\mathrm{loc}}(\mathbb R^3) \quad\text{for every }1\le p<2_\alpha^*.
		\end{gathered}
		\]
		The lower bound on the integral over \(B_R(0)\) implies \(u_0\ne0\). Set
		\[
		\xi_n=v_n-u_0.
		\]
		
		We give the nonlinear splittings used below. For each \(D\in\{A,m,B,L,Q\}\),
		\begin{equation}\label{eq4.12}
			D(v_n)=D(u_0)+D(\xi_n)+o_n(1).
		\end{equation}
		The assertions for \(A\) and \(m\) follow from weak convergence; the one for \(Q\) is the Brezis-Lieb identity. For the convex part of the logarithm, \eqref{eq2.3} and the \(\Delta_2\)-condition imply that, for every \(\eta>0\),
		\begin{equation}\label{eq4.13}
			|F_1(a+b)-F_1(a)| \le\eta F_1(a)+C_\eta F_1(b) \qquad\forall a,b\in\mathbb R.
		\end{equation}
		To verify this, first consider \(|b|\le\theta|a|\), where the mean value theorem and \(sF_1'(s)\le2F_1(s)\) give the first term when \(\theta\) is sufficiently small. On the complementary set, \(|a|\le|b|/\theta\), and the dilation estimate for \(F_1\) gives the second term.
		
		Apply \eqref{eq4.13} with \(a=\xi_n\), \(b=u_0\). The positive part of
		\[
		|F_1(v_n)-F_1(\xi_n)-F_1(u_0)|-\eta F_1(\xi_n)
		\]
		is dominated by \((C_\eta+1)F_1(u_0)\) and tends to zero almost everywhere. Since \(\int F_1(\xi_n)\) is bounded, dominated convergence gives
		\[
		\limsup_{n\to\infty} \int_{\mathbb R^3}|F_1(v_n)-F_1(\xi_n)-F_1(u_0)|\,dx \le C\eta.
		\]
		Letting \(\eta\downarrow0\) proves the splitting of \(\int F_1\). The bound \(|F_2'(s)|\le C|s|^{p-1}\), with \(2<p<q\), gives
		\[
		|F_2(a+b)-F_2(a)|\le\eta|a|^p+C_\eta|b|^p.
		\]
		The same argument proves the splitting of \(\int F_2\). Since \(L=2\int(F_2-F_1)\), the assertion for \(L\) follows.
		
		For the Poisson term, put \(p_\alpha=12/(3+2\alpha)\). The Hardy-Littlewood-Sobolev inequality gives
		\[
		\left|\int_{\mathbb R^3}(I_{2\alpha}*f)g\,dx\right| \le C\|f\|_{p_\alpha/2}\|g\|_{p_\alpha/2}.
		\]
		Local compactness and approximation of \(u_0\) in \(L^{p_\alpha}(\mathbb R^3)\) imply
		\[
		u_0\xi_n\to0\quad\text{in }L^{p_\alpha/2}(\mathbb R^3).
		\]
		Also \(\xi_n^2\rightharpoonup0\) in \(L^{p_\alpha/2}(\mathbb R^3)\), because this sequence is bounded and converges almost everywhere to zero. Expanding
		\[
		v_n^2=u_0^2+2u_0\xi_n+\xi_n^2
		\]
		in the continuous bilinear form above, all terms containing \(u_0\xi_n\) tend to zero, and
		\[
		\int_{\mathbb R^3}(I_{2\alpha}*u_0^2)\xi_n^2\,dx\to0.
		\]
		This proves the assertion for \(B\), and completes the proof of \eqref{eq4.12}. In particular,
		\begin{equation}\label{eq4.14}
			\begin{aligned}
				\mathcal G_0(v_n) &=\mathcal G_0(u_0)+\mathcal G_0(\xi_n)+o_n(1),\\
				\mathcal K_0(v_n) &=\mathcal K_0(u_0)+\mathcal K_0(\xi_n)+o_n(1).
			\end{aligned}
		\end{equation}
		
		We claim that \(\mathcal G_0(u_0)\le0\). Otherwise, \(\mathcal G_0(\xi_n)<0\) for large \(n\), since \(\mathcal G_0(v_n)=0\). By Lemma~\ref{Lem4.4}, there is \(s_n\in(0,1)\) such that \(T_{s_n}\xi_n\in\mathcal P_0\). Using \eqref{eq4.10}, we obtain
		\[
		\begin{aligned}
			c_0 &\le J_0(T_{s_n}\xi_n) =\mathcal K_0(T_{s_n}\xi_n)\\
			&\le\mathcal K_0(\xi_n) =c_0-\mathcal K_0(u_0)+o_n(1).
		\end{aligned}
		\]
		This contradicts \(u_0\ne0\) and \eqref{eq4.8}. Therefore \(\mathcal G_0(u_0)\le0\).
		
		Let \(t_0\in(0,1]\) be the unique parameter with \(T_{t_0}u_0\in\mathcal P_0\). The second splitting in \eqref{eq4.14} and \(\mathcal K_0(\xi_n)\ge0\) give
		\[
		\begin{aligned}
			c_0 \le\mathcal K_0(T_{t_0}u_0) \le\mathcal K_0(u_0)\le\liminf_{n\to\infty}\mathcal K_0(v_n)=c_0.
		\end{aligned}
		\]
		All inequalities are equalities. Strict monotonicity in \eqref{eq4.10} forces \(t_0=1\), so
		\[
		u_0\in\mathcal P_0, \qquad J_0(u_0)=c_0, \qquad\mathcal K_0(\xi_n)\to0.
		\]
		By \eqref{eq4.8},
		\[
		\xi_n\to0\quad\text{in }H^\alpha(\mathbb R^3)\cap L^q(\mathbb R^3).
		\]
		Writing \eqref{eq4.7} in terms of \(F_1,F_2\), and dropping its nonnegative quadratic and power terms, gives
		\[
		0\le\frac{2(3-\alpha)}{9-2\alpha} \int_{\mathbb R^3}F_1(\xi_n)\,dx \le\mathcal K_0(\xi_n)+C\int_{\mathbb R^3}F_2(\xi_n)\,dx\to0.
		\]
		Equations \eqref{eq2.2}-\eqref{eq2.3} imply convergence in \(L^{F_1}(\mathbb R^3)\). Consequently,
		\[
		v_n\to u_0\qquad\text{in }X.
		\]
		
		It remains to show that a minimizer is an unconstrained critical point. Let \(U\in\mathcal P_0\) satisfy \(J_0(U)=c_0\). The path \(t\mapsto T_tU\), including its endpoint at zero, has its unique maximum \(c_0\) at \(t=1\). Suppose that \(J_0'(U)\ne0\). Choose a small ball about \(U\), not containing zero, on which \(\|J_0'\|\) has a positive lower bound. A local descending deformation, constructed from a locally Lipschitz pseudogradient in the Banach space \(X\), lowers the energy on a smaller ball while leaving points outside the first ball fixed. On the compact part of the scaling path outside this smaller ball, the maximum is strictly below \(c_0\), by the uniqueness of the maximum. Choosing the deformation size below this gap produces an admissible path whose maximum is less than \(c_0\), contradicting \eqref{eq4.11}. Hence \(J_0'(U)=0\).
		
		Finally, if \(U\) changed sign, the fractional quadratic form would give
		\[
		A(|U|)<A(U),
		\]
		whereas \(m,B,L,Q\) are unchanged by taking the absolute value. Let \(t_*>0\) maximize \(J_0(T_t|U|)\). Then
		\[
		c_0\le J_0(T_{t_*}|U|) <J_0(T_{t_*}U) \le J_0(U)=c_0,
		\]
		a contradiction. Thus every minimizer has a fixed sign.
	\end{proof}
	
	\begin{lemma}\label{Lem4.6}
		The nonnegative minimizers in Proposition~\ref{Prop4.1} are strictly positive. After a choice of centers they form a compact set in $X$ and in $C^k$ on each bounded set, for every finite $k$. Their $L^\infty$ norms are uniformly bounded and they vanish uniformly at infinity after this centering.
	\end{lemma}
	\begin{proof}
		The resolvent argument in the proof of Lemma~\ref{Lem3.5} applies to the autonomous equation and gives uniform \(L^\infty\) bounds on bounded families of nonnegative solutions. Splitting the Riesz integral at distance one also gives
		\[
		\|\phi_U^\alpha\|_\infty\le C(\|U\|_\infty^2+\|U\|_2^2).
		\]
		The equation now has bounded right-hand side. On a fixed ball, subtract the Riesz potential of a cutoff of this right-hand side; the remainder is fractional harmonic on a smaller ball. The potential and interior harmonic estimates give uniform local $C^\sigma$ bounds for every nonintegral $\sigma<2\alpha$. Higher local regularity follows by repeating this decomposition with a H\"older right-hand side. Continuity suffices for positivity. If $U(x_0)=0$, use the zero function as an interior touching test in a small ball about $x_0$, and retain $U$ outside that ball. Its fractional Laplacian at $x_0$ is strictly negative, while the equation has right-hand side zero there. This contradicts the viscosity supersolution inequality. The passage from a continuous weak solution to this inequality is justified by the weak-viscosity equivalence \cite{ServadeiValdinoci}; it can be applied on any bounded ball.
		
		With the kernel \(K_\alpha\) defined in the proof of Lemma~\ref{Lem3.5}, the uniform \(L^\infty\) bound also gives $U\le C K_\alpha*U$. If centered solutions converge in $L^2$, the right-hand side converges uniformly, since $K_\alpha\in L^2$. For a fixed $L^2$ function its convolution with $K_\alpha$ is in $C_0$. This proves uniform vanishing and shows that maximum points cannot escape after the translations in Proposition~\ref{Prop4.1}. We may consequently first center a maximum at zero and obtain a compact set in $X$ and in $C^0$ on bounded sets. Strict positivity and compactness give a positive lower bound on each fixed ball. On such a ball the logarithm is smooth. Interior estimates for the coupled equations then bootstrap to uniform $C^k$ bounds. Compactness in $C^k$ on smaller balls follows by one additional regularity step.
	\end{proof}
	
	\begin{lemma}\label{Lem4.7}
		Every finite-energy weak solution \(u\in X\) of \eqref{eq4.2} satisfies
		\[
		\mathcal G_0(u)=0.
		\]
	\end{lemma}
	
	\begin{proof}
		We first justify the compactly supported test functions used in the identity. Applying the convex approximation of the absolute value in the fractional quadratic form gives the subsolution inequality for \(|u|\) used in Lemma~\ref{Lem3.5}. The same resolvent iteration therefore gives \(u\in L^\infty\). The Riesz representation gives \(\phi_u^\alpha\in L^\infty\), and interior estimates first yield \(u\in C^\beta_{\mathrm{loc}}\) for every nonintegral \(\beta<2\alpha\). In particular, one may choose \(\beta>1\). On bounded intervals the function \(s\mapsto s\log s^2\) is \(\sigma\)-H\"older for every \(\sigma\in(0,1)\). A second application of interior estimates to the coupled equations gives
		\[
		u\in C^{2\alpha+\sigma}_{\mathrm{loc}}(\mathbb R^3)
		\]
		for nonintegral \(2\alpha+\sigma\). Choose \(\sigma>1-\alpha\). Then the compactly supported first spatial variations of \(u\) belong to \(X\).
		
		Fix \(\chi\in C_c^\infty(\mathbb R^3)\), equal to one on \(B_1\) and zero outside \(B_2\). Set
		\[
		Y_R(x)=\chi(x/R)x,
		\]
		and let \(\Phi_s^R\) be its flow. With \(j_s^R=\det D\Phi_s^R\), define
		\[
		S_s^R u(x)=j_s^R(x)u(\Phi_s^R(x)).
		\]
		This deformation agrees with \(u\) outside \(B_{2R}\), and it is differentiable in \(X\) at zero. Its derivative is
		\[
		z_R(u)=(\operatorname{div}Y_R)u+Y_R\cdot\nabla u.
		\]
		Testing the weak equation with \(z_R(u)\) gives
		\begin{equation}\label{eq4.15}
			\left.\frac d{ds}J_0(S_s^R u)\right|_{s=0} =J_0'(u)[z_R(u)]=0.
		\end{equation}
		
		We pass to the limit as \(R\to\infty\) in this identity. The bounds
		\[
		\|DY_R\|_\infty\le C, \qquad \|D^2Y_R\|_\infty\le C/R
		\]
		are independent of \(R\). In the fractional quadratic form, make the change of variables \(\xi=\Phi_s^R(x)\), \(\zeta=\Phi_s^R(y)\). The derivatives at zero of the distorted kernel are bounded by \(C|\xi-\zeta|^{-3-2\alpha}\). The nonconstant Jacobian factors also produce differences bounded by
		\[
		C\min\{1,|\xi-\zeta|/R\}.
		\]
		Cauchy-Schwarz and the same radial integral as in Lemma~\ref{Lem3.3} show that their contribution is bounded by
		\[
		C R^{-\alpha}A(u)^{1/2}\|u\|_2 +C R^{-2\alpha}\|u\|_2^2.
		\]
		This tends to zero. The remaining integrand is dominated by a constant times the Gagliardo integrand and converges pointwise to the one for the full dilation. Thus
		\[
		\left.\frac d{ds}A(S_s^Ru)\right|_{s=0} \longrightarrow(3+2\alpha)A(u).
		\]
		
		For the Poisson term, the derivatives of the transformed Riesz kernel are bounded by a constant times \(|\xi-\zeta|^{-(3-2\alpha)}\). Dominated convergence gives
		\[
		\left.\frac d{ds}B(S_s^Ru)\right|_{s=0} \longrightarrow(9-2\alpha)B(u).
		\]
		For the mass and power terms the corresponding limits are \(3m(u)\) and \(3(q-1)Q(u)\). Finally, after the same change of variables in the logarithmic term, its first derivative is dominated by
		\[
		C\bigl(u^2+|u^2\log u^2|\bigr),
		\]
		which is integrable because \(u\in X\). Its limit is \(3L(u)+6m(u)\). Taking the limit in \eqref{eq4.15} now gives exactly \eqref{eq4.5}, and proves \(\mathcal G_0(u)=0\).
	\end{proof}
	
	\begin{theorem}\label{Thm4.1}
		Problem \eqref{eq4.2} has a positive ground state solution \(u_0\in X\) such that
		\[
		J_0(u_0)=c_0.
		\]
	\end{theorem}
	
	\begin{proof}
		By Proposition~\ref{Prop4.1}, there is a minimizer \(u_0\in\mathcal P_0\) satisfying
		\[
		J_0(u_0)=c_0, \qquad J_0'(u_0)=0.
		\]
		The same proposition shows that \(u_0\) has a fixed sign. Since \(J_0\) is even, we choose the sign so that \(u_0\ge0\). Lemma~\ref{Lem4.6} then gives
		\[
		u_0>0\qquad\text{in }\mathbb R^3.
		\]
		Let \(v\in X\setminus\{0\}\) be any other weak solution of \eqref{eq4.2}. Lemma~\ref{Lem4.7} gives \(\mathcal G_0(v)=0\), hence \(v\in\mathcal P_0\) and
		\[
		J_0(v)\ge c_0=J_0(u_0).
		\]
		Thus \(u_0\) is a positive ground state solution.
	\end{proof}
	
	\begin{lemma}\label{Lem4.8}
		For \(\mu\ge V_0\), define
		\[
		\begin{aligned}
			J_\mu(u) &=\frac12\int_{\mathbb R^3} \Bigl(|(-\Delta)^{\alpha/2}u|^2+(\mu+1)|u|^2\Bigr)\,dx +\frac14\int_{\mathbb R^3}\phi_u^\alpha u^2\,dx\\
			&\quad+\int_{\mathbb R^3}F_1(u)\,dx -\int_{\mathbb R^3}F_2(u)\,dx -\frac1q\int_{\mathbb R^3}|u|^q\,dx,
		\end{aligned}
		\]
		and
		\[
		\mathcal G_\mu(u)=\left.\frac d{dt}J_\mu(T_tu)\right|_{t=1}, \qquad \mathcal P_\mu=\{u\in X\setminus\{0\}:\mathcal G_\mu(u)=0\}, \qquad c_\mu=\inf_{\mathcal P_\mu}J_\mu.
		\]
		Thus \(J_{V_0}=J_0\) and \(c_{V_0}=c_0\). The map \(\mu\mapsto c_\mu\) is strictly increasing on \([V_0,\infty)\).
	\end{lemma}
	
	\begin{proof}
		The preceding autonomous arguments apply with \(V_0\) replaced by any \(\mu\ge V_0\). Let \(V_0\le\mu_1<\mu_2\), and take a ground state \(U\in\mathcal P_{\mu_2}\) with \(J_{\mu_2}(U)=c_{\mu_2}\). Let \(t>0\) be the unique parameter such that \(T_tU\in\mathcal P_{\mu_1}\). Then
		\[
		\begin{aligned}
			c_{\mu_1} &\le J_{\mu_1}(T_tU)\\
			&=J_{\mu_2}(T_tU) -\frac{\mu_2-\mu_1}{2}t^3\|U\|_2^2\\
			&<J_{\mu_2}(T_tU) \le\max_{s>0}J_{\mu_2}(T_sU) =J_{\mu_2}(U)=c_{\mu_2}.
		\end{aligned}
		\]
		This proves the assertion.
	\end{proof}
	
	\subsection{A local natural constraint}
	
	In this subsection, an integral without a specified domain is taken over \(\mathbb R^3\). We first define a smooth center on a neighborhood of the positive autonomous ground states and their translates. Take \(\theta\in C_c^\infty(\mathbb R^3)\), nonnegative and nonzero, and put
	\[
	r_u=\theta*u^2.
	\]
	For a compact centered family as in Lemma~\ref{Lem4.6}, \(r_U\) tends to zero uniformly at infinity and \(\inf_U\|r_U\|_\infty>0\). Choose \(\gamma>0\) smaller than half this infimum. Compactness supplies \(R_c>0\) such that
	\[
	r_U(x)<\gamma/2\qquad\text{for }|x|\ge R_c,
	\]
	uniformly in the centered family. The estimate
	\[
	\|r_u-r_v\|_\infty \le C\|u-v\|_2(\|u\|_2+\|v\|_2)
	\]
	extends these properties to a small \(X\) neighborhood of all their translates. In that neighborhood the set \(\{r_u>\gamma\}\) is contained in a ball of fixed radius about the nearby translate.
	
	Choose \(\vartheta\in C^\infty(\mathbb R)\), zero on \(( -\infty,\gamma]\) and positive on \((\gamma,\infty)\), and define
	\begin{equation}\label{eq4.16}
		\mathfrak b(u) =\frac{\displaystyle\int x \vartheta(r_u(x))\,dx} {\displaystyle\int \vartheta(r_u(x))\,dx}.
	\end{equation}
	The denominator has a uniform positive lower bound on a smaller neighborhood. Indeed, this holds first on the compact family, where \(\max r_U>2\gamma\), and then follows by uniform continuity. Differentiating gives
	\[
	D\mathfrak b(u)[v] =\frac{2\displaystyle\int(x-\mathfrak b(u))\vartheta'(r_u(x)) (\theta*(uv))(x)\,dx} {\displaystyle\int \vartheta(r_u(x))\,dx}.
	\]
	The variable \(x-\mathfrak b(u)\) is bounded on the support of the integrand, and
	\[
	\|\theta*(uv)\|_\infty\le C\|u\|_2\|v\|_2, \qquad \|\theta*(vw)\|_\infty\le C\|v\|_2\|w\|_2.
	\]
	These estimates and the analogous formula for the second derivative show that \(\mathfrak b\) is \(C^2\), with uniformly bounded first and second derivatives. It satisfies
	\[
	\mathfrak b(u(\cdot-z))=\mathfrak b(u)+z.
	\]
	Recenter the positive ground states by \(\mathfrak b\), and denote the resulting compact set by \(\mathcal S\). Thus \(\mathfrak b(U)=0\) for \(U\in\mathcal S\), and every positive ground state is a translate of a member of \(\mathcal S\).
	
	\begin{lemma}\label{Lem4.9}
		There is a $C^1$, translation equivariant vector field $Z$ on an $X$ neighborhood of $\{U(\cdot-z):U\in\mathcal S,\ z\in\mathbb R^3\}$, of the form
		\[
		Z(u)=\psi_u u,\qquad \operatorname{supp}\psi_u\subset B_{2R}(\mathfrak b(u)),
		\]
		where $R$ is fixed and $u\mapsto\psi_u$ is $C^1$ into $W^{1,\infty}$. On a smaller neighborhood the field and its first derivative are bounded. The function $\Gamma_0(u)=J_0'(u)[Z(u)]$ is $C^1$, and for some $\kappa>0$,
		\begin{equation}\label{eq4.17}
			\Gamma_0'(u)[Z(u)]\le-4\kappa, \qquad \mathcal G_0'(u)[Z(u)]\le-4\kappa.
		\end{equation}
	\end{lemma}
	\begin{proof}
		We justify first the negative directions at the minimizers. Let $Y_R(x)=\chi(x/R)x$, where $\chi=1$ on $B_1$, $\chi=0$ outside $B_2$, and let $\Phi_s^R$ be the flow of $Y_R$. Writing $j_s^R=\det D\Phi_s^R$, define the compactly supported deformation
		\[
		S_s^R U(x)=j_s^R(x)U(\Phi_s^R(x)).
		\]
		Outside $B_{2R}$ it is the identity. By Lemma~\ref{Lem4.6}, this path is twice differentiable in $X$ for fixed $R$, and its first derivative is
		\[
		z_R(U)=(\operatorname{div}Y_R)U+Y_R\cdot\nabla U.
		\]
		Since $J_0'(U)=0$, its second energy derivative is the quadratic form
		\begin{equation}\label{eq4.18}
			\begin{split}
				Q_U(v,v)={}&A(v)+V_0\|v\|_2^2+\int\phi_U^{\alpha}v^2\,dx +2\int(I_{2\alpha}*(Uv))Uv\,dx\\
				&-\int(\log U^2+2)v^2\,dx-(q-1)\int U^{q-2}v^2\,dx
			\end{split}
		\end{equation}
		evaluated at $v=z_R(U)$. This quadratic form is used only for compactly supported directions, on whose support $U$ has a positive lower bound.
		
		We next pass to the limit in the energy derivatives. We use the change-of-variables argument in Lemma~\ref{Lem4.7}, now keeping both time derivatives. The bounds
		\[
		\|DY_R\|_\infty\le C, \quad \|D^2Y_R\|_\infty\le C/R, \quad \|D^3Y_R\|_\infty\le C/R^2
		\]
		imply uniform bounds on the first two time derivatives of the Jacobians. The flow and its inverse are uniformly bi-Lipschitz for small \(|s|\). For \(\lambda=3+2\alpha\) and \(\lambda=3-2\alpha\), differentiating the distorted kernel gives
		\[
		\left|\left.\frac{d^j}{ds^j} | (\Phi_s^R)^{-1}\xi-(\Phi_s^R)^{-1}\zeta|^{-\lambda} \right|_{s=0}\right| \le C|\xi-\zeta|^{-\lambda},\qquad j=1,2.
		\]
		Indeed, the first and second time derivatives of the difference of the inverse flows are bounded by \(C|\xi-\zeta|\), using the bounds on \(DY_R\) and \(D^2Y_R\,Y_R\). For the kinetic term change variables $\xi=\Phi_s^R(x)$ and $\zeta=\Phi_s^R(y)$. Its numerator, after the Jacobians are divided out, is
		\[
		\left| \sqrt{\frac{j_s^R((\Phi_s^R)^{-1}\xi)} {j_s^R((\Phi_s^R)^{-1}\zeta)}}U(\xi) - \sqrt{\frac{j_s^R((\Phi_s^R)^{-1}\zeta)} {j_s^R((\Phi_s^R)^{-1}\xi)}}U(\zeta) \right|^2.
		\]
		The derivatives of the nonconstant Jacobian ratios are bounded by $C\min\{1,|\xi-\zeta|/R\}$. For such terms, put
		\[
		D_R(U)=\iint\frac{U(\zeta)^2 \min\{1,|\xi-\zeta|/R\}^2}{|\xi-\zeta|^{3+2\alpha}} \,d\xi\,d\zeta.
		\]
		The radial integration used in Lemma~\ref{Lem3.3} gives \(D_R(U)\le CR^{-2\alpha}\|U\|_2^2\). Cauchy-Schwarz bounds the mixed terms by \(C A(U)^{1/2}D_R(U)^{1/2}\). Thus the sum of all terms containing a Jacobian difference is bounded by
		\[
		CR^{-\alpha}A(U)^{1/2}\|U\|_2+CR^{-2\alpha}\|U\|_2^2.
		\]
		The remaining terms are dominated by a constant times the Gagliardo integrand, and their coefficients converge pointwise to those for the full dilation. For the Poisson term the dominating integrand is its positive Riesz integrand. For the other terms the majorants are $U^2$, $|U^2\log U^2|$ and $|U|^q$. Dominated convergence applies to each of the first two energy derivatives. We now make the uniformity on \(\mathcal S\) explicit. For the quadratic and Poisson terms, the differentiated expressions are bounded bilinear and quartic forms, with bounds independent of \(R\ge1\). Their differences on the bounded set \(\mathcal S\) are therefore bounded by \(C\|U-W\|_X\). The mass and power terms have the same bound, by the Sobolev embeddings. For the logarithmic term, \eqref{eq4.13} and the \(F_2\) growth bound imply that \(U\mapsto U^2\log U^2\) is continuous from \(X\) to \(L^1\). Its differentiated coefficients are uniformly bounded in \(R\). These estimates also apply to the first derivative of \(\mathcal G_0(S_s^R U)\).
		
		Let \(E_R(U)\) be the sum of the absolute errors in the two limits in \eqref{eq4.19}. The preceding estimates, including those for the full dilation, give
		\[
		|E_R(U)-E_R(W)| \le C\|U-W\|_X +C\|U^2\log U^2-W^2\log W^2\|_1, \qquad U,W\in\mathcal S,
		\]
		with \(C\) independent of \(R\ge1\). By compactness of \(\mathcal S\), the right-hand side tends to zero uniformly as \(\|U-W\|_X\to0\). Thus there is a modulus \(\omega\), independent of \(R\), such that
		\[
		|E_R(U)-E_R(W)|\le\omega(\|U-W\|_X), \qquad \omega(\rho)\to0\quad\text{as }\rho\downarrow0.
		\]
		For a fixed \(\rho>0\), take a finite \(\rho\)-net \(U_1,\ldots,U_N\) of \(\mathcal S\) in \(X\). Then
		\[
		\sup_{U\in\mathcal S}E_R(U) \le\max_{1\le j\le N}E_R(U_j)+\omega(\rho).
		\]
		The maximum tends to zero for fixed \(\rho\), by the pointwise limits already proved. Letting first \(R\to\infty\) and then \(\rho\downarrow0\) proves that the convergence is uniform on \(\mathcal S\). Hence
		\begin{equation}\label{eq4.19}
			Q_U(z_R(U),z_R(U))\longrightarrow \left.\frac{d^2}{ds^2}J_0(T_{e^s}U)\right|_{s=0}, \qquad \mathcal G_0'(U)[z_R(U)]\longrightarrow \left.\frac{d^2}{ds^2}J_0(T_{e^s}U)\right|_{s=0}.
		\end{equation}
		By \eqref{eq4.6} and $\mathcal G_0(U)=0$, the common limit is at most
		\[
		-\frac{(3+2\alpha)(6-4\alpha)}2 A(U)-3\alpha\|U\|_2^2.
		\]
		It is bounded above by a fixed negative number on $\mathcal S$. Fix $R$ large enough to keep both quantities in \eqref{eq4.19} negative.
		
		To extend these directions to arbitrary nearby functions, fix a smooth mollifier $\rho_\tau$. On $B_{2R}$ the functions in $\mathcal S$ have a uniform positive lower bound. For $\tau>0$ small set
		\begin{equation}\label{eq4.20}
			\begin{split}
				\psi_u(x)={}&\operatorname{div}Y_R(x-\mathfrak b(u))\\
				&+Y_R(x-\mathfrak b(u))\cdot \frac{\nabla(\rho_\tau*u)(x)}{(\rho_\tau*u)(x)}, \qquad Z(u)=\psi_u u.
			\end{split}
		\end{equation}
		In \eqref{eq4.20} the quotient is only used on the support of $Y_R$; there its denominator is bounded below on a sufficiently small neighborhood. Outside this range the denominator may be replaced by a fixed smooth positive extension without changing the formula on the support of $Y_R$. For fixed \(R\), Lemma~\ref{Lem4.6} gives uniform \(C^k\) compactness and a positive lower bound on \(B_{2R+1}\). Hence, as \(\tau\downarrow0\), the coefficients in \eqref{eq4.20} converge in \(W^{1,\infty}\) on this ball, and \(Z(U)\to z_R(U)\) in \(H^\alpha\) and in every finite \(L^p\), with fixed compact support, uniformly on \(\mathcal S\). The logarithmic coefficient in \eqref{eq4.18} is bounded on the support, and the Riesz form is continuous. Consequently both \(Q_U(Z(U),Z(U))\) and \(\mathcal G_0'(U)[Z(U)]\) converge uniformly to the corresponding expressions with \(z_R(U)\). Thus $\tau$ can be fixed so that both inequalities in \eqref{eq4.17} hold at the minimizers with a strict margin.
		
		Convolution and \eqref{eq4.16} show that $u\mapsto\psi_u$ is $C^1$ into $W^{1,\infty}$, with bounded derivative. Multiplication by such functions is bounded on $X$, proving the assertions about $Z$. To verify $\Gamma_0\in C^1$, expand it as
		\begin{equation}\label{eq4.21}
			\begin{split}
				\Gamma_0(u)={}& \int(-\Delta)^{\alpha/2}u(-\Delta)^{\alpha/2}(\psi_u u)\,dx +V_0\int\psi_u u^2\,dx+\int\phi_u^{\alpha}\psi_u u^2\,dx\\
				&-\int\psi_u u^2\log u^2\,dx-\int\psi_u|u|^q\,dx.
			\end{split}
		\end{equation}
		For the logarithmic term in \eqref{eq4.21}, write it as \(2\int\psi_u(F_2-F_1)(u)\). Its derivative in a direction \(v\in X\) is
		\[
		2\int D\psi_u[v](F_2-F_1)(u) +2\int\psi_u(F_2'-F_1')(u)v.
		\]
		The first integral is controlled by the \(W^{1,\infty}\) bound on \(D\psi_u\), and the second by the Orlicz H\"older inequality and \((f_2)\). These estimates also prove continuity of the derivative. The kinetic, mass, power and Poisson terms are treated by the bounded multiplier and Riesz estimates. Thus \(\Gamma_0\in C^1\) without requiring that \(J_0\) be twice differentiable on all of \(X\). At a critical point $U$, $\Gamma_0'(U)[Z(U)]=Q_U(Z(U),Z(U))$, since the term $J_0'(U)[DZ(U)Z(U)]$ vanishes. Continuity, compactness and translation equivariance extend the strict bounds to a uniform neighborhood.
	\end{proof}
	
	We now transfer this construction to \(J_\varepsilon\). Fix \(\delta>0\) so small that \(M_\delta\subset\Lambda\). Choose \(d>0\) with \(4d<\delta\). Its size will be reduced below using the continuity of \(V\) near \(M\). Let \(\chi_d\in C_c^\infty(\mathbb R^3)\) satisfy
	\[
	0\le\chi_d\le1, \qquad \chi_d=1\text{ on }B_d, \qquad \chi_d=0\text{ outside }B_{2d}.
	\]
	For \(y\in M_{2d}\) and \(U\in\mathcal S\), define
	\begin{equation}\label{eq4.22}
		W_{\varepsilon,y,U}(x) =\chi_d(\varepsilon x-y)U(x-y/\varepsilon).
	\end{equation}
	These functions are supported in \(\Lambda_\varepsilon\). The compactness of \(\mathcal S\), Lemma~\ref{Lem2.2}, and \(F_1(\chi_d u)\le F_1(u)\) imply
	\begin{equation}\label{eq4.23}
		\sup_{y\in M_{2d},\,U\in\mathcal S} \|W_{\varepsilon,y,U}(\cdot+y/\varepsilon)-U\|_X\to0.
	\end{equation}
	For the Orlicz part, this follows first for a fixed \(U\) by dominated convergence and then uniformly on \(\mathcal S\) by a finite covering and the \(\Delta_2\)-condition. The same finite-covering argument applies to the Sobolev cutoff estimate.
	
	Since \(V\) is bounded on the compact set \(M_{4d}\),
	\[
	\sup_{0<\varepsilon<1,\,y\in M_{2d},\,U\in\mathcal S} \|W_{\varepsilon,y,U}\|_\varepsilon<\infty.
	\]
	For fixed \(\varepsilon>0\), the family in \eqref{eq4.22} is compact in \(X_\varepsilon\): the parameter sets are compact, translation and cutoff multiplication are continuous in \(X\), and the weighted term is controlled on a common bounded support. Put
	\[
	d_\varepsilon(u) =\inf_{y\in M_{2d},\,U\in\mathcal S} \|u-W_{\varepsilon,y,U}\|_\varepsilon.
	\]
	This is a Lipschitz function on \(X_\varepsilon\). By \eqref{eq4.23}, a fixed small neighborhood of this family is contained in the neighborhood of translates on which the field \(Z\) and the center \(\mathfrak b\) are defined.
	
	For \(u\) in this neighborhood, set
	\[
	\Gamma_\varepsilon(u)=J_\varepsilon'(u)[Z(u)].
	\]
	If \(\operatorname{dist}(\varepsilon\mathfrak b(u),M)<2d\) and \(\varepsilon\) is small, then
	\[
	\operatorname{supp}\psi_u \subset B_{2R}(\mathfrak b(u))\subset\Lambda_\varepsilon.
	\]
	The original and penalized nonlinearities consequently agree on this support. Thus
	\begin{equation}\label{eq4.24}
		\Gamma_\varepsilon(u)-\Gamma_0(u) =\int_{\mathbb R^3}(V(\varepsilon x)-V_0)\psi_u u^2\,dx.
	\end{equation}
	This identity also establishes \(\Gamma_\varepsilon\in C^1\) on the neighborhood. Its derivative in the direction \(v\) is the corresponding derivative of \(\Gamma_0\) plus
	\[
	\int_{\mathbb R^3}(V(\varepsilon x)-V_0) \bigl(D\psi_u[v]u^2+2\psi_u uv\bigr)\,dx.
	\]
	No derivative of \(V\) occurs. The support and multiplier bounds give
	\[
	|\Gamma_\varepsilon(u)-\Gamma_0(u)| +\|\Gamma_\varepsilon'(u)-\Gamma_0'(u)\|_{X_\varepsilon'} \le C\sup_{x\in M_{3d}}|V(x)-V_0|.
	\]
	Because \(V=V_0\) on the compact set \(M\), the supremum tends to zero as \(d\downarrow0\). We first choose \(d\), then an outer neighborhood radius, and finally \(\varepsilon_0>0\), so that Lemma~\ref{Lem4.9} implies, throughout this outer neighborhood,
	\begin{equation}\label{eq4.25}
		\Gamma_\varepsilon'(u)[Z(u)]\le-2\kappa, \qquad \mathcal G_0'(u)[Z(u)]\le-2\kappa \qquad(0<\varepsilon<\varepsilon_0).
	\end{equation}
	The functions \(\Gamma_\varepsilon'\), \(Z\), and \(DZ\) are bounded there, uniformly in \(\varepsilon\), with the norms of \(X_\varepsilon\) and its dual.
	
	Let \(\mathcal F_su\) denote the local flow of \(Z\), so that
	\[
	\frac d{ds}\mathcal F_su=Z(\mathcal F_su), \qquad\mathcal F_0u=u.
	\]
	The uniform bounds give a time interval \([-s_0,s_0]\) on a smaller neighborhood, independent of small \(\varepsilon\), and
	\begin{equation}\label{eq4.26}
		\|\mathcal F_su-u\|_\varepsilon\le C|s|.
	\end{equation}
	Since \(Z(v)=\psi_vv\), the integral equation also gives, almost everywhere,
	\[
	\mathcal F_su(x) =u(x)\exp\left(\int_0^s\psi_{\mathcal F_\tau u}(x)\,d\tau\right).
	\]
	In particular, this flow preserves the sign and support of a function.
	
	Choose an inner radius \(r>0\) small enough that its closure lies in the smaller neighborhood and
	\[
	|\mathcal G_0(u)|<\kappa s_0
	\]
	there. This is possible because \(\mathcal G_0\) vanishes on the limiting ground-state family and \eqref{eq4.23} is uniform. Decrease \(\varepsilon_0\) if needed so that the physical center moves by less than \(d\) along the flow; indeed,
	\[
	|\varepsilon\mathfrak b(\mathcal F_su) -\varepsilon\mathfrak b(u)|\le C\varepsilon|s|.
	\]
	Define
	\begin{equation}\label{eq4.27}
		\begin{aligned}
			\mathcal U_\varepsilon &=\{u\in X_\varepsilon:d_\varepsilon(u)<r, \ \operatorname{dist}(\varepsilon\mathfrak b(u),M)<d\},\\
			\Sigma_\varepsilon &=\{u\in\mathcal U_\varepsilon:\Gamma_\varepsilon(u)=0\}.
		\end{aligned}
	\end{equation}
	The radius is also chosen so that this neighborhood is disjoint from zero. Compactness of \(\mathcal S\) and its disjointness from zero allow this choice. By \eqref{eq4.25}, \(\Sigma_\varepsilon\) is a \(C^1\) hypersurface of \(X_\varepsilon\).
	
	\begin{proposition}\label{Prop4.2}
		For all sufficiently small \(\varepsilon>0\), the following hold. Every critical point of \(J_\varepsilon|_{\Sigma_\varepsilon}\) is an unconstrained critical point of \(J_\varepsilon\), and
		\begin{equation}\label{eq4.28}
			\|J_\varepsilon'(u)\|_{X_\varepsilon'} \le C\|(J_\varepsilon|_{\Sigma_\varepsilon})'(u)\| \qquad\forall u\in\Sigma_\varepsilon.
		\end{equation}
		For every \(u\in\Sigma_\varepsilon\), there is a unique \(s(u)\in(-s_0,s_0)\) such that
		\[
		p(u):=\mathcal F_{s(u)}u\in\mathcal P_0.
		\]
		Moreover,
		\begin{equation}\label{eq4.29}
			c_0\le J_0(p(u))\le J_\varepsilon(p(u))\le J_\varepsilon(u),
		\end{equation}
		and
		\begin{equation}\label{eq4.30}
			J_\varepsilon(u)-J_\varepsilon(p(u))\ge\kappa|s(u)|^2.
		\end{equation}
	\end{proposition}
	
	\begin{proof}
		Let \(u\) be a constrained critical point. Since \(\Gamma_\varepsilon'(u)\ne0\), the Lagrange multiplier rule gives
		\[
		J_\varepsilon'(u)=\lambda\Gamma_\varepsilon'(u).
		\]
		Acting on \(Z(u)\), and using \(\Gamma_\varepsilon(u)=0\), we find
		\[
		0=\lambda\Gamma_\varepsilon'(u)[Z(u)].
		\]
		The first inequality in \eqref{eq4.25} forces \(\lambda=0\). Thus \(J_\varepsilon'(u)=0\).
		
		More generally, given \(v\in X_\varepsilon\), define
		\[
		v_T=v- \frac{\Gamma_\varepsilon'(u)[v]} {\Gamma_\varepsilon'(u)[Z(u)]}Z(u).
		\]
		Then \(\Gamma_\varepsilon'(u)[v_T]=0\), so \(v_T\in T_u\Sigma_\varepsilon\). The uniform bounds and \eqref{eq4.25} give
		\[
		\|v_T\|_\varepsilon\le C\|v\|_\varepsilon.
		\]
		Since \(J_\varepsilon'(u)[Z(u)]=0\),
		\[
		J_\varepsilon'(u)[v]=J_\varepsilon'(u)[v_T].
		\]
		Taking the supremum over \(\|v\|_\varepsilon\le1\) proves \eqref{eq4.28}. In particular, a constrained Palais-Smale sequence is a Palais-Smale sequence for the original functional.
		
		Now fix \(u\in\Sigma_\varepsilon\). By the second inequality in \eqref{eq4.25},
		\[
		\frac d{ds}\mathcal G_0(\mathcal F_su) \le-2\kappa.
		\]
		The choice \(|\mathcal G_0(u)|<\kappa s_0\) gives
		\[
		\mathcal G_0(\mathcal F_{-s_0}u)>0, \qquad \mathcal G_0(\mathcal F_{s_0}u)<0.
		\]
		Thus there is exactly one \(s(u)\in(-s_0,s_0)\) with \(\mathcal G_0(\mathcal F_{s(u)}u)=0\). The projected function is nonzero, because the flow preserves zero sets and \(u\ne0\).
		
		Put \(e(s)=J_\varepsilon(\mathcal F_su)\). The identity \(e'(s)=\Gamma_\varepsilon(\mathcal F_su)\) and \(\Gamma_\varepsilon\in C^1\) show that \(e\in C^2\), even though only \(J_\varepsilon\in C^1\) is used. We have
		\[
		e'(0)=0, \qquad e''(s)=\Gamma_\varepsilon'(\mathcal F_su) [Z(\mathcal F_su)]\le-2\kappa.
		\]
		Integrating twice, for either sign of \(s\), gives
		\[
		e(s)\le e(0)-\kappa s^2.
		\]
		Taking \(s=s(u)\) proves \eqref{eq4.30}. Finally, \(p(u)\in\mathcal P_0\) gives \(J_0(p(u))\ge c_0\), and \eqref{eq4.3} gives \(J_\varepsilon(p(u))\ge J_0(p(u))\). This proves \eqref{eq4.29}.
	\end{proof}
	
	\subsection{Asymptotic behavior}
	
	We next compare functions of energy close to \(c_0\) on \(\Sigma_\varepsilon\) with the autonomous ground states. The comparison also keeps the relevant sublevels away from the boundary of \(\mathcal U_\varepsilon\).
	
	\begin{proposition}\label{Prop4.3}
		Let \(\varepsilon_n\to0\), and let \(u_n\in\Sigma_{\varepsilon_n}\) satisfy
		\[
		J_{\varepsilon_n}(u_n)\to c_0.
		\]
		Then
		\begin{equation}\label{eq4.31}
			d_{\varepsilon_n}(u_n)\to0, \qquad \operatorname{dist}(\varepsilon_n\mathfrak b(u_n),M)\to0.
		\end{equation}
		Moreover, there exist translations \(y_n\in\mathbb R^3\), a point \(y_0\in M\), and a positive ground state \(U\in\mathcal S\) such that, up to a subsequence,
		\[
		u_n(\cdot+y_n)\to U\quad\text{in }X, \qquad\varepsilon_ny_n\to y_0.
		\]
		There also exist \(\eta>0\) and \(\varepsilon_1>0\) such that
		\begin{equation}\label{eq4.32}
			\begin{aligned}
				\{u\in\Sigma_\varepsilon:J_\varepsilon(u)\le c_0+2\eta\} \subset \{u:d_\varepsilon(u)<r/2,\ \operatorname{dist}(\varepsilon\mathfrak b(u),M)<d/2\}
			\end{aligned}
		\end{equation}
		for \(0<\varepsilon<\varepsilon_1\). On these sublevels, \(J_\varepsilon|_{\Sigma_\varepsilon}\) satisfies the Palais-Smale condition.
	\end{proposition}
	
	\begin{proof}
		Let
		\[
		s_n=s(u_n),\qquad p_n=p(u_n)=\mathcal F_{s_n}u_n.
		\]
		By \eqref{eq4.29},
		\[
		c_0\le J_0(p_n)\le J_{\varepsilon_n}(p_n) \le J_{\varepsilon_n}(u_n)=c_0+o_n(1).
		\]
		Consequently,
		\begin{equation}\label{eq4.33}
			J_0(p_n)\to c_0, \qquad J_{\varepsilon_n}(p_n)-J_0(p_n)\to0.
		\end{equation}
		Also, by \eqref{eq4.30},
		\[
		0\le\kappa|s_n|^2 \le J_{\varepsilon_n}(u_n)-J_{\varepsilon_n}(p_n)\to0.
		\]
		Thus \(s_n\to0\), and \eqref{eq4.26} gives
		\begin{equation}\label{eq4.34}
			\|u_n-p_n\|_{\varepsilon_n}\le C|s_n|\to0.
		\end{equation}
		
		Since \(p_n\in\mathcal P_0\), Proposition~\ref{Prop4.1} applies. There are translations \(y_n\) and a fixed-sign autonomous minimizer \(U\) such that
		\[
		p_n(\cdot+y_n)\to U\qquad\text{in }X.
		\]
		The limit cannot be negative. To see this, choose nonnegative \(W_n'\) from the family \eqref{eq4.22} with \(\|u_n-W_n'\|_{\varepsilon_n}<r\). Pointwise, \(p_n^-\le|p_n-W_n'|\), hence
		\[
		\|p_n^-\|_2\le C r+o_n(1).
		\]
		On the other hand, if \(U<0\), strong \(L^2\) convergence would give \(\|p_n^-\|_2\to\|U\|_2\). The \(L^2\) norms of the autonomous minimizers have a positive infimum, by their compactness up to translations and their nonzero limits. Choosing \(r\) below this fixed bound excludes a negative \(U\). Recenter \(U\) so that \(U\in\mathcal S\).
		
		By the continuity and translation equivariance of \(\mathfrak b\),
		\[
		\mathfrak b(p_n)-y_n\to0.
		\]
		The bounded derivative of \(\mathfrak b\) and \eqref{eq4.34} also give
		\[
		\mathfrak b(u_n)-\mathfrak b(p_n)\to0.
		\]
		Since \(\varepsilon_n\mathfrak b(u_n)\in M_d\), the sequence \(\varepsilon_ny_n\) is bounded. Passing to a subsequence,
		\[
		\varepsilon_ny_n\to y_0\in M_d.
		\]
		
		The two nonnegative terms in \eqref{eq4.3} and \eqref{eq4.33} imply
		\begin{equation}\label{eq4.35}
			\int_{\mathbb R^3}(V(\varepsilon_nx)-V_0)p_n^2\,dx\to0.
		\end{equation}
		Fix \(R>0\). After translating by \(y_n\), the potential converges uniformly on \(B_R(0)\):
		\[
		V(\varepsilon_nx+\varepsilon_ny_n)\to V(y_0).
		\]
		Strong \(L^2\) convergence on this ball and \eqref{eq4.35} yield
		\[
		0\le (V(y_0)-V_0)\int_{B_R(0)}U^2\,dx\le0.
		\]
		As \(U\ne0\), we obtain \(V(y_0)=V_0\). Since \(M_d\subset\Lambda\), this means \(y_0\in M\). Thus the center assertion in \eqref{eq4.31} holds.
		
		It remains to prove closeness in the weighted norm, rather than only in \(X\). Put \(z_n=\varepsilon_ny_n\). For large \(n\), \(z_n\in M_{2d}\), so
		\[
		W_n=W_{\varepsilon_n,z_n,U}
		\]
		belongs to the family defining \(d_{\varepsilon_n}\). By \eqref{eq4.23},
		\[
		p_n-W_n\to0\qquad\text{in }X.
		\]
		For the weighted term, use
		\[
		\begin{aligned}
			&\int_{\mathbb R^3}(V(\varepsilon_nx)-V_0)(p_n-W_n)^2\,dx\\
			&\qquad\le2\int_{\mathbb R^3}(V(\varepsilon_nx)-V_0)p_n^2\,dx +2\int_{\mathbb R^3}(V(\varepsilon_nx)-V_0)W_n^2\,dx.
		\end{aligned}
		\]
		The first term tends to zero by \eqref{eq4.35}. In the second, change variables \(x\mapsto x+y_n\):
		\[
		\int_{\mathbb R^3} \bigl(V(\varepsilon_nx+z_n)-V_0\bigr) \chi_d(\varepsilon_nx)^2U(x)^2\,dx.
		\]
		On a fixed ball, this tends to zero by \(z_n\to y_0\in M\). Outside the ball it is bounded by \(C\int_{B_R^c}U^2\), because \(V\) is bounded on \(M_{4d}\), which contains the physical support of the integrand. First let \(n\to\infty\), then \(R\to\infty\). The second term therefore also tends to zero. We have proved
		\[
		\|p_n-W_n\|_{\varepsilon_n}\to0.
		\]
		Combining this with \eqref{eq4.34},
		\[
		d_{\varepsilon_n}(u_n) \le\|u_n-p_n\|_{\varepsilon_n} +\|p_n-W_n\|_{\varepsilon_n}\to0.
		\]
		This proves \eqref{eq4.31} and the asserted strong convergence after translation.
		
		Suppose that \eqref{eq4.32} failed for every choice of \(\eta,\varepsilon_1>0\). We could then choose \(\varepsilon_n<1/n\) and \(u_n\in\Sigma_{\varepsilon_n}\) with
		\[
		J_{\varepsilon_n}(u_n)\le c_0+1/n,
		\]
		such that either \(d_{\varepsilon_n}(u_n)\ge r/2\) or \(\operatorname{dist}(\varepsilon_n\mathfrak b(u_n),M)\ge d/2\). By \eqref{eq4.29}, their energies converge to \(c_0\), contradicting \eqref{eq4.31}. Thus fixed \(\eta,\varepsilon_1>0\) exist.
		
		Finally, fix \(\varepsilon\in(0,\varepsilon_1)\), and let \((v_n)\subset\Sigma_\varepsilon\) be a constrained Palais-Smale sequence in one of these sublevels. It is bounded in \(X_\varepsilon\) by the definition of the neighborhood. Equation \eqref{eq4.28} gives \(J_\varepsilon'(v_n)\to0\). Lemma~\ref{Lem3.4} then gives, after a subsequence,
		\[
		v_n\to v\qquad\text{in }X_\varepsilon.
		\]
		By \eqref{eq4.32}, its limit satisfies \(d_\varepsilon(v)\le r/2<r\) and \(\operatorname{dist}(\varepsilon\mathfrak b(v),M)\le d/2<d\). It therefore lies in \(\mathcal U_\varepsilon\), where continuity gives \(\Gamma_\varepsilon(v)=0\). Hence \(v\in\Sigma_\varepsilon\), which proves the constrained Palais-Smale condition.
	\end{proof}
	
	\begin{theorem}\label{Thm4.2}
		Let \(\varepsilon_n\to0\), and let \(u_n\in\Sigma_{\varepsilon_n}\) be a sequence of nonnegative critical points such that
		\[
		J_{\varepsilon_n}(u_n)\to c_0, \qquad J_{\varepsilon_n}'(u_n)=0.
		\]
		Then there exist a sequence \(\{y_n\}\subset\mathbb R^3\) and a positive ground state solution \(w\) of \eqref{eq4.2} such that, up to a subsequence,
		\begin{itemize}
			\item[(i)]
			\[
			w_n(x)=u_n(x+y_n)\to w\qquad\text{in }X;
			\]
			\item[(ii)] there exists \(y_0\in M\) such that
			\[
			\varepsilon_ny_n\to y_0;
			\]
			\item[(iii)] there exists \(C>0\) such that
			\[
			\|w_n\|_{L^\infty(\mathbb R^3)}\le C\qquad\forall n,
			\]
			the convergence in (i) is locally uniform, and
			\[
			w_n(x)\to0\quad\text{uniformly in }n\text{ as }|x|\to\infty.
			\]
		\end{itemize}
	\end{theorem}
	
	\begin{proof}
		Proposition~\ref{Prop4.3} gives translations \(y_n\), \(y_0\in M\), and a positive ground state \(w\) such that
		\[
		w_n=u_n(\cdot+y_n)\to w\quad\text{in }X, \qquad\varepsilon_ny_n\to y_0.
		\]
		This proves (i) and (ii). In particular, for every \(p\in[2,2_\alpha^*)\),
		\[
		w_n\to w\qquad\text{in }L^p(\mathbb R^3).
		\]
		The strong convergence also gives the tail estimates
		\[
		\lim_{R\to\infty}\limsup_{n\to\infty} \int_{B_R^c(0)}\left( |(-\Delta)^{\alpha/2}w_n|^2+|w_n|^2\right)\,dx=0,
		\]
		\[
		\lim_{R\to\infty}\limsup_{n\to\infty} \int_{B_R^c(0)}|w_n|^p\,dx=0,
		\]
		and
		\[
		\lim_{R\to\infty}\limsup_{n\to\infty} \int_{B_R^c(0)}F_1(w_n)\,dx=0.
		\]
		For the last limit, use \(F_1(w_n)\le C(F_1(w)+F_1(w_n-w))\), \(\int F_1(w_n-w)\to0\), and \(F_1(w)\in L^1\).
		
		The functions \(w_n\) solve the translated auxiliary equations. Since \(V\ge0\), \(\phi_{w_n}^\alpha\ge0\), and \eqref{eq2.4} is uniform in the spatial variable,
		\[
		\bigl(( -\Delta)^\alpha+1\bigr)w_n \le C w_n+C w_n^{q-1}.
		\]
		The neighborhood defining \(\Sigma_{\varepsilon_n}\) bounds \(u_n\) uniformly in \(H^\alpha\). The integrability iteration in Lemma~\ref{Lem3.5} therefore gives
		\[
		\sup_n\|w_n\|_\infty<\infty.
		\]
		It also gives a uniform bound for \(\phi_{w_n}^\alpha\).
		
		Fix \(R>0\). Because \(\varepsilon_ny_n\to y_0\in M\), for large \(n\) the set \(\varepsilon_n B_{2R}(0)+\varepsilon_ny_n\) lies in \(\Lambda\). Moreover,
		\[
		\sup_n\sup_{x\in B_{2R}(0)} V(\varepsilon_nx+\varepsilon_ny_n)<\infty.
		\]
		All right-hand sides are thus uniformly bounded on \(B_{2R}(0)\). Interior H\"older estimates, the global \(L^2\) bound, and local compactness give a subsequence converging uniformly on \(B_R(0)\). Its limit must be \(w\), by the strong \(L^2\) convergence. A diagonal argument proves local uniform convergence on \(\mathbb R^3\).
		
		We finally prove the uniform decay, retaining the nonlocal contribution from the core. Since \(w_n\) is uniformly bounded, \(w_n^{q-1}\le Cw_n\), and therefore
		\[
		\bigl(( -\Delta)^\alpha+1\bigr)w_n\le Cw_n.
		\]
		The resolvent comparison gives
		\[
		0\le w_n\le C K_\alpha*w_n.
		\]
		As \(K_\alpha\in L^2\) and \(w_n\to w\) in \(L^2\),
		\[
		\|K_\alpha*(w_n-w)\|_\infty \le\|K_\alpha\|_2\|w_n-w\|_2\to0.
		\]
		The convolution of two \(L^2\) functions is continuous and vanishes at infinity: approximate each factor in \(L^2\) by compactly supported smooth functions and use the same inequality. Consequently,
		\[
		\lim_{R\to\infty}\limsup_{n\to\infty} \sup_{|x|\ge R}w_n(x)=0.
		\]
		For each fixed \(n\), the comparison with \(K_\alpha*w_n\) also gives \(w_n(x)\to0\) as \(|x|\to\infty\). Including the finitely many initial indices yields the decay uniformly for the entire sequence, and proves (iii).
	\end{proof}
	
	\section{Multiplicity and concentration}\label{Sec5}
	
	\subsection{Multiplicity of solutions}
	
	In this subsection, we prove the existence of multiple solutions to \eqref{eq2.5} by means of Lusternik-Schnirelmann category theory. Set
	\[
	M:=\{x\in\Lambda:\ V(x)=V_{0}\}, \qquad M_{\delta}:=\{x\in\mathbb{R}^{3}:\ \operatorname{dist}(x,M)\le \delta\},
	\]
	where \(\delta>0\) is chosen so small that \(M_{\delta}\subset\Lambda\).
	
	We first recall the following Lusternik-Schnirelmann theorem in the form for compact subsets of a \(C^1\) Finsler manifold \cite[Theorem 3.1]{Szulkin}.
	
	\begin{theorem}\label{Thm5.1}
		Let \(\mathscr M\) be a \(C^1\) Finsler manifold, and let \(I\in C^1(\mathscr M,\mathbb R)\) be bounded from below. Assume that every sublevel set
		\[
		I^d:=\{u\in\mathscr M:I(u)\le d\}
		\]
		is complete in the induced metric. Let \(K\subset\mathscr M\) be compact, and suppose that \(I\) satisfies the \((PS)_c\)-condition for every \(c\le\sup_K I\). Then \(I\) has at least
		\[
		\operatorname{cat}_{\mathscr M}(K)
		\]
		distinct critical points with values not exceeding \(\sup_K I\).
	\end{theorem}
	
	Here the category is defined using closed subsets contractible in \(\mathscr M\). For completeness, the minimax values in this theorem are
	\[
	c_j=\inf_{A\in\mathscr A_j}\sup_{u\in A}I(u), \qquad \mathscr A_j=\{A\subset\mathscr M:A\text{ compact},\ \operatorname{cat}_{\mathscr M}(A)\ge j\}.
	\]
	For \(j\le\operatorname{cat}_{\mathscr M}(K)\), these classes are nonempty and \(c_j\le\sup_K I\). The conclusion follows from \cite[Theorem 3.1]{Szulkin}, which applies to these minimax classes of compact sets on a \(C^1\) Finsler manifold and also treats coincident minimax values.
	
	We now construct the standard family of test functions. Let \(u_0\) be the positive ground state solution of \eqref{eq4.2} given by Theorem~\ref{Thm4.1}, translated so that \(\mathfrak b(u_0)=0\). Fix \(\varphi\in C^\infty([0,\infty))\) such that
	\[
	0\le \varphi\le 1, \qquad \varphi(t)=1 \ \text{for } 0\le t\le \frac{\delta}{2}, \qquad \varphi(t)=0 \ \text{for } t\ge \delta.
	\]
	For \(y\in M\), define
	\[
	\Upsilon_{\varepsilon,y}(x):= \varphi(|\varepsilon x-y|)\, u_0\!\left(x-\frac{y}{\varepsilon}\right).
	\]
	These functions belong to \(X_\varepsilon\) and are supported in \(\Lambda_\varepsilon\). They approach the cutoff family in \eqref{eq4.22}, uniformly in \(y\), in the norm of \(X_\varepsilon\). Also
	\[
	\sup_{y\in M}|\Gamma_\varepsilon(\Upsilon_{\varepsilon,y})|\to0.
	\]
	By \eqref{eq4.25}, there is a unique small parameter \(a_\varepsilon(y)\) such that
	\[
	\Gamma_\varepsilon\bigl(\mathcal F_{a_\varepsilon(y)} \Upsilon_{\varepsilon,y}\bigr)=0.
	\]
	Indeed, along this flow the derivative of \(\Gamma_\varepsilon\) is at most \(-2\kappa\), and
	\[
	|a_\varepsilon(y)| \le \frac{|\Gamma_\varepsilon(\Upsilon_{\varepsilon,y})|}{2\kappa}.
	\]
	The parameter is continuous in \(y\), by strict monotonicity and continuous dependence of the flow. We set
	\[
	\Phi_\varepsilon(y):= \mathcal F_{a_\varepsilon(y)}\Upsilon_{\varepsilon,y} \in\Sigma_\varepsilon.
	\]
	For all small \(\varepsilon\), membership in the neighborhood defining \(\Sigma_\varepsilon\) follows from the uniform cutoff convergence, \eqref{eq4.26}, and continuity of \(\mathfrak b\).
	
	Choose \(\rho>0\) such that \(M_\delta\subset B_\rho(0)\) and define \(\zeta:\mathbb R^3\to\mathbb R^3\) by
	\[
	\zeta(x)=
	\begin{cases}
		x, & |x|\le \rho,\\[2pt]
		\rho\,\dfrac{x}{|x|}, & |x|\ge \rho.
	\end{cases}
	\]
	Fix \(p\in(2,2_\alpha^*)\), and define
	\[
	\beta_\varepsilon:X_\varepsilon\setminus\{0\}\to\mathbb R^3, \qquad \beta_\varepsilon(u):= \frac{\displaystyle\int_{\mathbb R^3}\zeta(\varepsilon x)|u(x)|^p\,dx} {\displaystyle\int_{\mathbb R^3}|u(x)|^p\,dx}.
	\]
	Clearly,
	\[
	\beta_\varepsilon(tu)=\beta_\varepsilon(u) \qquad \forall\, t>0,\ u\neq0.
	\]
	We distinguish this truncated barycenter from the local center \(\mathfrak b\) in \eqref{eq4.16}. On an \(L^p\)-bounded family with \(L^p\) norms bounded away from zero,
	\begin{equation}\label{eq5.1}
		|\beta_\varepsilon(u)-\beta_\varepsilon(v)| \le C\|u-v\|_p,
	\end{equation}
	uniformly in \(\varepsilon\). This follows by subtracting the two quotients, using \(|\zeta|\le\rho\), and applying
	\[
	\bigl\||u|^p-|v|^p\bigr\|_1 \le C_p(\|u\|_p^{p-1}+\|v\|_p^{p-1})\|u-v\|_p.
	\]
	
	\begin{lemma}\label{Lem5.1}
		The following limit holds:
		\[
		\lim_{\varepsilon\to0}J_\varepsilon(\Phi_\varepsilon(y))=c_0 \qquad \text{uniformly in } y\in M.
		\]
		In particular, if we set
		\[
		h(\varepsilon):=\varepsilon+ \sup_{y\in M}|J_\varepsilon(\Phi_\varepsilon(y))-c_0|,
		\]
		then \(h(\varepsilon)\to0\) and
		\[
		\Phi_\varepsilon(M)\subset \widetilde\Sigma_\varepsilon:= \{u\in\Sigma_\varepsilon:J_\varepsilon(u)\le c_0+h(\varepsilon)\}
		\]
		for \(\varepsilon>0\) sufficiently small.
	\end{lemma}
	
	\begin{proof}
		Let \(\varepsilon_n\to0\) and \(y_n\in M\). Setting
		\[
		\psi_n(z):=\varphi(\varepsilon_n|z|)\,u_0(z),
		\]
		we have, after the change of variables \(x=z+\frac{y_n}{\varepsilon_n}\),
		\[
		\Upsilon_{\varepsilon_n,y_n}(x)=\psi_n(z).
		\]
		Lemma~\ref{Lem2.2} gives convergence of the fractional Sobolev part. Since \(0\le\varphi\le1\), the convexity of \(F_1\) gives
		\[
		F_1((1-\varphi(\varepsilon_n|z|))u_0(z))\le F_1(u_0(z)).
		\]
		Dominated convergence and \eqref{eq2.3} yield
		\[
		\psi_n\to u_0\qquad\text{in }X.
		\]
		On the support of \(\psi_n\), one has \(\varepsilon_nz+y_n\in M_\delta\subset\Lambda\). Thus the penalized and exact nonlinearities coincide there.
		
		The potential term converges as well. For fixed \(R>0\), continuity of \(V\) and compactness of \(M\) give
		\[
		\sup_{y\in M,\,|z|\le R}|V(\varepsilon_nz+y)-V_0|\to0.
		\]
		On the complement of \(B_R\), where \(\psi_n\ne0\), the potential is bounded by its maximum on \(M_\delta\). Therefore
		\[
		\begin{aligned}
			0&\le\int_{\mathbb R^3} (V(\varepsilon_nz+y_n)-V_0)\psi_n^2\,dz\\
			&\le o_n(1)\|u_0\|_2^2 +C\int_{|z|>R}u_0^2\,dz.
		\end{aligned}
		\]
		First taking \(n\to\infty\), then \(R\to\infty\), shows that this integral tends to zero. The Poisson term converges by Lemma~\ref{Lem2.1}; the logarithmic terms converge in their \(F_1,F_2\) decomposition, and the power term converges in \(L^q\). Consequently,
		\[
		J_{\varepsilon_n}(\Upsilon_{\varepsilon_n,y_n})\to J_0(u_0)=c_0.
		\]
		The same cutoff convergence, the continuity of the field in \eqref{eq4.20}, and \eqref{eq4.24} imply
		\[
		\Gamma_{\varepsilon_n}(\Upsilon_{\varepsilon_n,y_n})\to0.
		\]
		Hence \(a_{\varepsilon_n}(y_n)\to0\), and \eqref{eq4.26} gives
		\[
		\|\Phi_{\varepsilon_n}(y_n) -\Upsilon_{\varepsilon_n,y_n}\|_{\varepsilon_n}\to0.
		\]
		Along this short flow, the energy derivative is \(\Gamma_{\varepsilon_n}\), which is uniformly bounded on the outer neighborhood. Integrating it, we obtain
		\[
		J_{\varepsilon_n}(\Phi_{\varepsilon_n}(y_n)) -J_{\varepsilon_n}(\Upsilon_{\varepsilon_n,y_n})\to0.
		\]
		This proves the asserted convergence. Since the sequences \(\varepsilon_n\to0\) and \(y_n\in M\) were arbitrary, it is uniform in \(y\). Finally, \eqref{eq4.29} gives \(J_\varepsilon(\Phi_\varepsilon(y))\ge c_0\), and the last statement follows from the definition of \(h\).
	\end{proof}
	
	\begin{lemma}\label{Lem5.2}
		Let \(\varepsilon_n\to0\) and let \(u_n\in\widetilde\Sigma_{\varepsilon_n}\). Then there exist \(v_n\in\mathcal P_0\cap X_{\varepsilon_n}\), translations \(y_n\in\mathbb R^3\), a point \(y_0\in M\), and a positive ground state solution \(w\) of \eqref{eq4.2} such that, up to a subsequence,
		\[
		\|u_n-v_n\|_{\varepsilon_n}\to0, \qquad J_0(v_n)\to c_0, \qquad J_{\varepsilon_n}(v_n)-J_0(v_n)\to0,
		\]
		\[
		u_n(\cdot+y_n)\to w, \qquad v_n(\cdot+y_n)\to w \qquad\text{in }X,
		\]
		and
		\[
		\varepsilon_n y_n\to y_0.
		\]
	\end{lemma}
	
	\begin{proof}
		Take the projection in Proposition~\ref{Prop4.2}, and put
		\[
		v_n=p(u_n)=\mathcal F_{s(u_n)}u_n\in\mathcal P_0.
		\]
		By \eqref{eq4.29} and Lemma~\ref{Lem5.1},
		\[
		\begin{aligned}
			c_0\le J_0(v_n) \le J_{\varepsilon_n}(v_n) \le J_{\varepsilon_n}(u_n) \le c_0+h(\varepsilon_n).
		\end{aligned}
		\]
		It follows that
		\[
		J_0(v_n)\to c_0, \qquad J_{\varepsilon_n}(v_n)-J_0(v_n)\to0.
		\]
		Moreover, \eqref{eq4.30} implies
		\[
		\kappa|s(u_n)|^2 \le J_{\varepsilon_n}(u_n)-J_{\varepsilon_n}(v_n) \le h(\varepsilon_n)\to0.
		\]
		The flow estimate therefore gives
		\[
		\|u_n-v_n\|_{\varepsilon_n}\to0.
		\]
		Since \((v_n)\) is a minimizing sequence on \(\mathcal P_0\), Proposition~\ref{Prop4.1} provides translations \(y_n\) and a nonzero autonomous minimizer \(w\) for which
		\[
		v_n(\cdot+y_n)\to w\qquad\text{in }X.
		\]
		The same convergence holds for \(u_n\). The choice of the neighborhood around the positive family excludes a negative limiting minimizer, as in Proposition~\ref{Prop4.3}; thus \(w>0\).
		
		For clarity, the location of the centers follows from the potential defect in \eqref{eq4.3}. It gives
		\[
		\int_{\mathbb R^3} (V(\varepsilon_nx)-V_0)v_n^2\,dx\to0.
		\]
		The local-center continuity gives \(\mathfrak b(v_n)-y_n\to\mathfrak b(w)\), while the definition of \(\mathcal U_{\varepsilon_n}\) and the short flow keep \(\varepsilon_n\mathfrak b(v_n)\) in a fixed neighborhood of \(M\). Consequently, \((\varepsilon_ny_n)\) is bounded. Along a subsequence, let \(\varepsilon_ny_n\to y_0\). For a ball on which \(\int w^2>0\), local strong convergence and continuity of \(V\) then give
		\[
		(V(y_0)-V_0)\int_{B_R}w^2\,dx=0.
		\]
		The point \(y_0\) lies inside \(\Lambda\), so \(V(y_0)=V_0\) implies \(y_0\in M\).
	\end{proof}
	
	For the rest of this subsection, choose \(\eta>0\) as in \eqref{eq4.32} and set
	\[
	\mathscr M_\varepsilon =\{u\in\Sigma_\varepsilon:J_\varepsilon(u)<c_0+\eta\}, \qquad f_\varepsilon(u) =\frac{J_\varepsilon(u)-c_0}{c_0+\eta-J_\varepsilon(u)}.
	\]
	Here \(\Sigma_\varepsilon\) is defined in \eqref{eq4.27}. The manifold is equipped with the Finsler norm induced by \(\|\cdot\|_\varepsilon\). Since \(J_\varepsilon\ge c_0\) on \(\Sigma_\varepsilon\), the functional \(f_\varepsilon\) is nonnegative and belongs to \(C^1(\mathscr M_\varepsilon,\mathbb R)\). For small \(\varepsilon\), \(\widetilde\Sigma_\varepsilon\subset\mathscr M_\varepsilon\).
	
	\begin{lemma}\label{Lem5.3}
		For every sufficiently small fixed \(\varepsilon>0\), each finite sublevel of \(f_\varepsilon\) is complete in the induced Finsler metric, and \(f_\varepsilon\) satisfies the \((PS)\)-condition on these sublevels.
	\end{lemma}
	
	\begin{proof}
		Fix \(A\ge0\). The sublevel \(f_\varepsilon\le A\) is equivalently
		\[
		J_\varepsilon(u)\le c_0+\frac{\eta A}{1+A}<c_0+\eta.
		\]
		By \eqref{eq4.32}, all its points satisfy
		\[
		d_\varepsilon(u)<\frac r2, \qquad \operatorname{dist}(\varepsilon\mathfrak b(u),M)<\frac d2.
		\]
		Let \((u_n)\) be Cauchy in the induced metric on this sublevel. The ambient norm distance is bounded by the length of any connecting curve, so \((u_n)\) is Cauchy in \(X_\varepsilon\). Let \(u_n\to u\) there. Continuity gives
		\[
		d_\varepsilon(u)\le\frac r2<r, \qquad \operatorname{dist}(\varepsilon\mathfrak b(u),M)\le\frac d2<d, \qquad \Gamma_\varepsilon(u)=0.
		\]
		Thus \(u\in\Sigma_\varepsilon\), and
		\[
		J_\varepsilon(u)\le c_0+\frac{\eta A}{1+A}.
		\]
		The point belongs to \(\mathscr M_\varepsilon\), not to its energy boundary. Here the sublevel is equipped with the restriction of the metric of \(\mathscr M_\varepsilon\). In a local \(C^1\) graph chart around \(u\), the Finsler norm is equivalent to the norm on the tangent hyperplane. Joining the chart coordinates by straight segments shows that ambient convergence implies convergence in the induced metric. This proves completeness of the sublevel.
		
		Now let \((u_n)\) be a Palais-Smale sequence for \(f_\varepsilon\) in this sublevel. On it,
		\[
		\frac1\eta\le \frac{\eta}{(c_0+\eta-J_\varepsilon(u_n))^2} \le\frac{(1+A)^2}{\eta}.
		\]
		Since
		\[
		f_\varepsilon'(u) =\frac{\eta}{(c_0+\eta-J_\varepsilon(u))^2} (J_\varepsilon|_{\Sigma_\varepsilon})'(u),
		\]
		we have
		\[
		\|(J_\varepsilon|_{\Sigma_\varepsilon})'(u_n)\|\to0.
		\]
		Equation \eqref{eq4.28} implies \(J_\varepsilon'(u_n)\to0\) in \(X_\varepsilon'\). Lemma~\ref{Lem3.4} then gives strong convergence in \(X_\varepsilon\), after taking a subsequence. The same boundary argument keeps the limit in \(\mathscr M_\varepsilon\), and local charts give convergence in its induced metric. Thus the required \((PS)\)-condition holds.
	\end{proof}
	
	\begin{lemma}\label{Lem5.4}
		The barycenter map satisfies
		\[
		\sup_{u\in\widetilde\Sigma_\varepsilon} \operatorname{dist}\bigl(\beta_\varepsilon(u),M\bigr)\to0 \qquad\text{as }\varepsilon\to0.
		\]
		The neighborhood radii can be chosen so that, for every sufficiently small \(\varepsilon>0\),
		\[
		\beta_\varepsilon(\mathscr M_\varepsilon)\subset M_\delta.
		\]
	\end{lemma}
	
	\begin{proof}
		Assume by contradiction that there exist \(\delta_0>0\), \(\varepsilon_n\to0\), and \(u_n\in\widetilde\Sigma_{\varepsilon_n}\) such that
		\[
		\operatorname{dist}\bigl(\beta_{\varepsilon_n}(u_n),M\bigr) \ge\delta_0\qquad\forall n.
		\]
		Applying Lemma~\ref{Lem5.2}, we find \(v_n\in\mathcal P_0\), a sequence \(\{y_n\}\subset\mathbb R^3\), a point \(y_0\in M\), and a positive ground state solution \(w\) of \eqref{eq4.2} such that
		\[
		\|u_n-v_n\|_{\varepsilon_n}\to0, \qquad \widetilde w_n(x):=v_n(x+y_n)\to w\quad\text{in }X, \qquad \varepsilon_n y_n\to y_0.
		\]
		The \(L^p\) norms of \(u_n\) and \(v_n\) stay bounded away from zero. Since \(u_n-v_n\to0\) in \(L^p(\mathbb R^3)\), \eqref{eq5.1} gives
		\[
		\beta_{\varepsilon_n}(u_n)-\beta_{\varepsilon_n}(v_n)\to0.
		\]
		By the change of variables \(x=\xi+y_n\),
		\[
		\beta_{\varepsilon_n}(v_n) = \frac{\displaystyle\int_{\mathbb R^3} \zeta(\varepsilon_n\xi+\varepsilon_n y_n)|\widetilde w_n(\xi)|^p\,d\xi} {\displaystyle\int_{\mathbb R^3}|\widetilde w_n(\xi)|^p\,d\xi}.
		\]
		The denominator tends to \(\|w\|_p^p>0\). For the numerator, write its difference from \(y_0\|\widetilde w_n\|_p^p\) as the sum of the integrals over \(B_R\) and \(B_R^c\). On \(B_R\),
		\[
		\sup_{|\xi|\le R} |\zeta(\varepsilon_n\xi+\varepsilon_n y_n)-y_0|\to0,
		\]
		because \(\zeta(y_0)=y_0\). On \(B_R^c\), the absolute value of the integral is bounded by
		\[
		2\rho\int_{|\xi|>R}|\widetilde w_n(\xi)|^p\,d\xi.
		\]
		Strong \(L^p\) convergence makes this uniformly small as \(R\to\infty\). Consequently,
		\[
		\beta_{\varepsilon_n}(v_n)\to y_0\in M.
		\]
		Hence \(\beta_{\varepsilon_n}(u_n)\to y_0\), a contradiction.
		
		We verify the range assertion on the whole manifold, which is needed for the homotopies below. Uniformly for \(y\in M_{2d}\) and \(U\in\mathcal S\), the same core and tail calculation gives
		\[
		\beta_\varepsilon(W_{\varepsilon,y,U})-y\to0.
		\]
		Here compactness of \(\mathcal S\) gives uniform \(L^p\) tails and a positive lower bound for the denominator. If \(u\in\mathcal U_\varepsilon\), choose \(W_{\varepsilon,y,U}\) with \(\|u-W_{\varepsilon,y,U}\|_\varepsilon<r\). For small \(r\), the same denominator lower bound holds for \(u\). By \eqref{eq5.1} and the Sobolev embedding,
		\[
		\operatorname{dist}(\beta_\varepsilon(u),M) \le Cr+2d+o_\varepsilon(1).
		\]
		We impose \(Cr<\delta/4\) when choosing the inner radius, before fixing \(\eta\) in \eqref{eq4.32}. Since \(4d<\delta\), the right-hand side is less than \(\delta\) for all sufficiently small \(\varepsilon\). This proves the range assertion for \(\mathscr M_\varepsilon\subset\mathcal U_\varepsilon\).
	\end{proof}
	
	\begin{lemma}\label{Lem5.5}
		The barycenter map satisfies
		\[
		\lim_{\varepsilon\to0}\beta_\varepsilon(\Phi_\varepsilon(y))=y \qquad \text{uniformly in } y\in M.
		\]
	\end{lemma}
	
	\begin{proof}
		Let \(\varepsilon_n\to0\) and \(y_n\in M\). The flow estimate in the proof of Lemma~\ref{Lem5.1} and \eqref{eq5.1} give
		\[
		\beta_{\varepsilon_n}(\Phi_{\varepsilon_n}(y_n)) - \beta_{\varepsilon_n}(\Upsilon_{\varepsilon_n,y_n})\to0.
		\]
		Using the change of variables \(x=z+\frac{y_n}{\varepsilon_n}\), we obtain
		\[
		\beta_{\varepsilon_n}(\Upsilon_{\varepsilon_n,y_n}) = \frac{\displaystyle\int_{\mathbb R^3}\zeta(\varepsilon_n z+y_n) |\varphi(\varepsilon_n|z|)u_0(z)|^p\,dz} {\displaystyle\int_{\mathbb R^3}|\varphi(\varepsilon_n|z|)u_0(z)|^p\,dz}.
		\]
		Since \(y_n\in M\subset B_\rho(0)\), we have \(\zeta(y_n)=y_n\). The denominator converges to \(\int_{\mathbb R^3}|u_0|^p\,dz>0\) by dominated convergence. For the difference between the numerator and \(y_n\) times the denominator, the integral over \(B_R\) is bounded by
		\[
		\varepsilon_nR\int_{B_R}|u_0|^p\,dz,
		\]
		since \(\zeta\) is Lipschitz with constant one. The complementary integral is bounded by
		\[
		2\rho\int_{|z|>R}|u_0|^p\,dz.
		\]
		Taking \(n\to\infty\), then \(R\to\infty\), yields
		\[
		\beta_{\varepsilon_n}(\Upsilon_{\varepsilon_n,y_n})-y_n\to0.
		\]
		It follows that \(\beta_{\varepsilon_n}(\Phi_{\varepsilon_n}(y_n))-y_n\to0\). Since the sequences \(\varepsilon_n\to0\) and \(y_n\in M\) were arbitrary, the convergence is uniform in \(y\in M\).
	\end{proof}
	
	\begin{lemma}\label{Lem5.6}
		Assume \((V_1)\)-\((V_2)\). Then there exists \(\varepsilon_1>0\) such that, for every \(\varepsilon\in(0,\varepsilon_1)\), the auxiliary problem \eqref{eq2.5} has at least
		\[
		\operatorname{cat}_{M_\delta}(M)
		\]
		distinct nontrivial critical points in \(\widetilde\Sigma_\varepsilon\).
	\end{lemma}
	
	\begin{proof}
		We apply Theorem~\ref{Thm5.1} to
		\[
		\mathscr M=\mathscr M_\varepsilon, \qquad I=f_\varepsilon, \qquad K=K_\varepsilon:=\Phi_\varepsilon(M).
		\]
		The set \(K_\varepsilon\) is compact, since \(M\) is compact and \(\Phi_\varepsilon\) is continuous. By Lemma~\ref{Lem5.1},
		\[
		K_\varepsilon\subset\widetilde\Sigma_\varepsilon \subset\mathscr M_\varepsilon
		\]
		for small \(\varepsilon\). The functional \(f_\varepsilon\) is nonnegative, and Lemma~\ref{Lem5.3} supplies the completeness of its sublevels and the required \((PS)\)-condition.
		
		It remains to prove that
		\[
		\operatorname{cat}_{\mathscr M_\varepsilon}(K_\varepsilon) \ge\operatorname{cat}_{M_\delta}(M).
		\]
		Let \(A_1,\ldots,A_n\) be closed subsets of \(\mathscr M_\varepsilon\) covering \(K_\varepsilon\), each contractible in \(\mathscr M_\varepsilon\). Let
		\[
		h_i:[0,1]\times A_i\to\mathscr M_\varepsilon
		\]
		be a contraction of \(A_i\) to a point \(z_i\in\mathscr M_\varepsilon\). Set
		\[
		B_i=\Phi_\varepsilon^{-1}(A_i)\subset M.
		\]
		These are closed subsets of \(M\), and
		\[
		M=\bigcup_{i=1}^n B_i.
		\]
		Define
		\[
		H_i(t,y)=
		\begin{cases}
			(1-2t)y+2t\,\beta_\varepsilon(\Phi_\varepsilon(y)), &0\le t\le \frac12,\\[1mm]
			\beta_\varepsilon\bigl(h_i(2t-1,\Phi_\varepsilon(y))\bigr), &\frac12\le t\le1.
		\end{cases}
		\]
		For the first half of the homotopy, Lemma~\ref{Lem5.5} gives
		\[
		|H_i(t,y)-y|\le |\beta_\varepsilon(\Phi_\varepsilon(y))-y|<\delta.
		\]
		For the second half, Lemma~\ref{Lem5.4} applies on the whole \(\mathscr M_\varepsilon\). Therefore
		\[
		H_i([0,1]\times B_i)\subset M_\delta,
		\]
		with
		\[
		H_i(0,y)=y, \qquad H_i(1,y)=\beta_\varepsilon(z_i).
		\]
		Thus each \(B_i\) is contractible in \(M_\delta\), and
		\[
		\operatorname{cat}_{M_\delta}(M)\le n.
		\]
		Taking the smallest such number of sets proves the category inequality.
		
		Theorem~\ref{Thm5.1} now gives at least \(\operatorname{cat}_{M_\delta}(M)\) distinct critical points of \(f_\varepsilon\) with values at most \(\sup_{K_\varepsilon}f_\varepsilon\). The transformation from \(J_\varepsilon\) to \(f_\varepsilon\) is strictly increasing, with positive derivative. These points are therefore critical points of \(J_\varepsilon|_{\Sigma_\varepsilon}\), and their energies satisfy
		\[
		c_0\le J_\varepsilon(u) \le\sup_{y\in M}J_\varepsilon(\Phi_\varepsilon(y)) \le c_0+h(\varepsilon).
		\]
		Proposition~\ref{Prop4.2} gives \(J_\varepsilon'(u)=0\). The points are nonzero because \(\mathcal U_\varepsilon\) is disjoint from zero. This completes the proof.
	\end{proof}
	
	\begin{lemma}\label{Lem5.7}
		The inner neighborhood radius can be chosen so that every critical point of \(J_\varepsilon\) in \(\Sigma_\varepsilon\) is positive, for every sufficiently small \(\varepsilon>0\).
	\end{lemma}
	
	\begin{proof}
		Write \(u^-=\max\{-u,0\}\). Testing the equation with \(-u^-\), using the positivity of the Poisson term and the convexity of \(F_1\), gives
		\[
		\|u^-\|_{H_\varepsilon}^2 \le \int_{\mathbb R^3}g(\varepsilon x,u^-)u^-\,dx.
		\]
		For the fractional part, the sign is justified by the pointwise inequality
		\[
		(a-b)(-a^-+b^-)\ge |a^- - b^-|^2.
		\]
		By \eqref{eq2.4}, for a fixed \(p\in(2,q)\),
		\[
		\left(1-\frac{b_0}{V_0+1}\right) \|u^-\|_{H_\varepsilon}^2 \le C\|u^-\|_p^p+\|u^-\|_q^q.
		\]
		The Sobolev embeddings, uniformly in \(\varepsilon\), imply
		\[
		C\|u^-\|_p^p+\|u^-\|_q^q \le C\bigl(\|u^-\|_p^{p-2}+\|u^-\|_q^{q-2}\bigr) \|u^-\|_{H_\varepsilon}^2.
		\]
		Choose a nonnegative cutoff ground state \(W_{\varepsilon,y,U}\) with \(\|u-W_{\varepsilon,y,U}\|_\varepsilon<r\). Pointwise,
		\[
		u^-\le |u-W_{\varepsilon,y,U}|.
		\]
		Hence
		\[
		\|u^-\|_p+\|u^-\|_q\le Cr.
		\]
		In choosing \(r\), require also that
		\[
		C\bigl((Cr)^{p-2}+(Cr)^{q-2}\bigr) <1-\frac{b_0}{V_0+1}.
		\]
		The preceding inequalities then force \(u^-=0\). Since the neighborhood excludes zero, \(u\) is a nontrivial nonnegative critical point. Lemma~\ref{Lem3.5} yields \(u>0\) in \(\mathbb R^3\).
	\end{proof}
	
	\subsection{Proof of Theorem~\ref{Thm1.1}}
	
	\begin{proof}
		By Lemma~\ref{Lem5.6}, there exists \(\varepsilon_1>0\) such that, for every \(\varepsilon\in(0,\varepsilon_1)\), the functional \(J_\varepsilon\) has at least
		\[
		\operatorname{cat}_{M_\delta}(M)
		\]
		distinct nontrivial critical points in \(\widetilde\Sigma_\varepsilon\). By Lemma~\ref{Lem5.7}, these critical points are positive for all sufficiently small \(\varepsilon\). Choose \(\varepsilon_*>0\) so that both conclusions hold and \(h(\varepsilon)\to0\) as in Lemma~\ref{Lem5.1}.
		
		Fix \(\varepsilon\in(0,\varepsilon_*)\), and denote the resulting critical points by
		\[
		w_{\varepsilon,1},\ldots,w_{\varepsilon,k}, \qquad k\ge\operatorname{cat}_{M_\delta}(M).
		\]
		They satisfy
		\[
		J_\varepsilon'(w_{\varepsilon,j})=0, \qquad c_0\le J_\varepsilon(w_{\varepsilon,j})\le c_0+h(\varepsilon), \qquad w_{\varepsilon,j}>0.
		\]
		
		We next prove that
		\begin{equation}\label{eq5.2}
			w_{\varepsilon,j}(x)<a_0 \qquad \forall x\in \mathbb R^3\setminus \Lambda_\varepsilon, \quad j=1,\dots,k.
		\end{equation}
		Assume by contradiction that there exist \(\varepsilon_n\to0\), points \(x_n\in\mathbb R^3\setminus\Lambda_{\varepsilon_n}\), and positive critical points \(\widetilde w_n\) chosen among the above family such that
		\[
		\widetilde w_n(x_n)\ge a_0.
		\]
		Since
		\[
		c_0\le J_{\varepsilon_n}(\widetilde w_n)\le c_0+h(\varepsilon_n),
		\]
		we have
		\[
		J_{\varepsilon_n}(\widetilde w_n)\to c_0.
		\]
		Applying Theorem~\ref{Thm4.2}, we find a sequence \(\{y_n\}\subset\mathbb R^3\), a point \(y_0\in M\), and a positive ground state solution \(w\) of \eqref{eq4.2} such that, up to a subsequence,
		\[
		\widehat w_n(x):=\widetilde w_n(x+y_n)\to w \qquad \text{in }X,
		\]
		\[
		\varepsilon_n y_n\to y_0,
		\]
		and
		\[
		\widehat w_n(x)\to0 \qquad \text{uniformly in }n \text{ as } |x|\to\infty.
		\]
		Since \(y_0\in M\subset\Lambda\), there exists \(r_*>0\) such that
		\[
		B_{2r_*}(y_0)\subset \Lambda.
		\]
		Thus, for \(n\) large,
		\[
		B_{r_*/\varepsilon_n}(y_n)\subset \Lambda_{\varepsilon_n}.
		\]
		The uniform decay of \(\widehat w_n\) then gives
		\[
		\widetilde w_n(x)<a_0 \qquad \forall x\in\mathbb R^3\setminus B_{r_*/\varepsilon_n}(y_n),
		\]
		for \(n\) large. Since
		\[
		\mathbb R^3\setminus \Lambda_{\varepsilon_n} \subset \mathbb R^3\setminus B_{r_*/\varepsilon_n}(y_n),
		\]
		this contradicts the choice of \(x_n\). Hence \eqref{eq5.2} holds.
		
		Consequently, each \(w_{\varepsilon,j}\) is a positive solution of the original problem \eqref{eq1.3}. Defining
		\[
		v_{\varepsilon,j}(x):=w_{\varepsilon,j}\!\left(\frac{x}{\varepsilon}\right), \qquad j=1,\dots,k,
		\]
		and \(\phi_{\varepsilon,j}(x)=\phi_{w_{\varepsilon,j}}^\alpha(x/\varepsilon)\), we conclude that problem \eqref{eq1.2} admits at least
		\[
		\operatorname{cat}_{M_\delta}(M)
		\]
		distinct positive solutions for every \(\varepsilon\in(0,\varepsilon_*)\).
		
		Each of these solutions is continuous and tends to zero at infinity, by the kernel comparison in Theorem~\ref{Thm4.2}, so it attains its positive maximum. It remains to prove the concentration of the maximum points. Let \(v_\varepsilon\) be one of the above positive solutions of \eqref{eq1.2}, and let \(\eta_\varepsilon\in\mathbb R^3\) be a global maximum point of \(v_\varepsilon\). Set
		\[
		w_\varepsilon(x):=v_\varepsilon(\varepsilon x).
		\]
		Then \(w_\varepsilon\) is one of the positive critical points of \(J_\varepsilon\), and
		\[
		x_\varepsilon:=\frac{\eta_\varepsilon}{\varepsilon}
		\]
		is a global maximum point of \(w_\varepsilon\).
		
		Now let \(\varepsilon_n\to0\). Since
		\[
		c_0\le J_{\varepsilon_n}(w_{\varepsilon_n})\le c_0+h(\varepsilon_n),
		\]
		we have
		\[
		J_{\varepsilon_n}(w_{\varepsilon_n})\to c_0.
		\]
		Applying Theorem~\ref{Thm4.2}, we obtain, up to a subsequence, a sequence \(\{y_n\}\subset\mathbb R^3\), a point \(y_0\in M\), and a positive ground state solution \(w\) of \eqref{eq4.2} such that
		\[
		\widetilde w_n(x):=w_{\varepsilon_n}(x+y_n)\to w \qquad \text{in }X,
		\]
		\[
		\varepsilon_n y_n\to y_0,
		\]
		and
		\[
		\widetilde w_n(x)\to0 \qquad \text{uniformly in }n \text{ as } |x|\to\infty.
		\]
		
		Set \(x_n=x_{\varepsilon_n}\). We claim that \(\{x_n-y_n\}\) is bounded. Suppose otherwise. Since \(w>0\), there exist \(R>0\) and \(a>0\) such that
		\[
		\sup_{x\in B_R(0)}w(x)\ge 2a.
		\]
		By Theorem~\ref{Thm4.2}(iii), \(\{\widetilde w_n\}\) is uniformly bounded in \(L^\infty(\mathbb R^3)\). Standard interior Hölder estimates then imply, up to a subsequence,
		\[
		\widetilde w_n\to w \qquad \text{uniformly on }B_R(0).
		\]
		Hence, for \(n\) large,
		\[
		\sup_{x\in B_R(0)}\widetilde w_n(x)\ge a.
		\]
		Since \(x_n\) is a global maximum point of \(w_{\varepsilon_n}\),
		\[
		w_{\varepsilon_n}(x_n)\ge \sup_{x\in B_R(0)}\widetilde w_n(x)\ge a.
		\]
		But
		\[
		w_{\varepsilon_n}(x_n)=\widetilde w_n(x_n-y_n),
		\]
		and \(|x_n-y_n|\to+\infty\), which contradicts the uniform decay of \(\widetilde w_n\). Thus \(\{x_n-y_n\}\) is bounded.
		
		It follows that
		\[
		\eta_{\varepsilon_n} = \varepsilon_n x_n = \varepsilon_n y_n+\varepsilon_n(x_n-y_n)\to y_0\in M.
		\]
		Hence
		\[
		V(\eta_{\varepsilon_n})\to V(y_0)=V_0.
		\]
		Since \(\varepsilon_n\to0\) was arbitrary, we conclude that
		\[
		\operatorname{dist}(\eta_\varepsilon,M)\to0, \qquad \lim_{\varepsilon\to0}V(\eta_\varepsilon)=V_0.
		\]
		This completes the proof.
	\end{proof}
	
	\section*{Acknowledgment}
	The authors thank the anonymous referee for the careful reading and helpful comments. This work was supported by the National Natural Science Foundation of China (12301145, 12561020, 12261107) and the Yunnan Fundamental Research Projects (202401AU070123, 202601AT070048).
	
	{\bf Data availability:} Data sharing is not applicable to this article as no new data were created or analyzed in this study.
	
	{\bf Conflict of Interests:} The authors declare that there is no conflict of interest.

\end{document}